\documentclass[12pt]{amsart}

\newcommand{\bff}{\bf}

\newcounter{sigmasigma}
\usepackage{graphics}
\usepackage{amssymb}
\usepackage{amsxtra}

\input diagrams

\usepackage[arrow, matrix]{xy}

\usepackage{ifthen}

\theoremstyle{plain}

\newtheorem{theo}{Theorem}[section]
\newtheorem{lemm}[theo]{Lemma}
\newtheorem{prop}[theo]{Proposition}
\newtheorem{coro}[theo]{Corollary}

\theoremstyle{definition}

\newtheorem{defi}[theo]{Definition}
\newtheorem{rema}[theo]{Remark}

\newtheorem{exam}[theo]{Example}

\newfont{\rmm}{cmr10 scaled 1000}
\newfont{\itt}{cmsl10 scaled 1000}

\newfont{\rM}{cmr10 scaled 1700}

\newcommand{\lb}{\label}
\newcommand{\mlb}{\label}

\newcommand{\rrf}[1]{(\ref{#1})}
\newcommand{\mrf}[1]{\ref{#1}}

\newenvironment{lis}
{
\begin{list}
{\textbf{\rm \bf\arabic{sigmasigma}.}}
{\setlength{\leftmargin}{.3in}
\usecounter{sigmasigma}}
}
{
\setcounter{sigmasigma}{0}
\end{list}
}

\newcommand{\beli}{\begin{lis}}
\def\enli{\end{lis}}

\newenvironment{liss}
{
\begin{list}
{\textbf{\rm \arabic{sigmasigma})}}
{\setlength{\leftmargin}{.3in}
\usecounter{sigmasigma}}
}
{
\setcounter{sigmasigma}{0}
\end{list}
}

\newcommand{\belis}{\begin{liss}}
\def\enlis{\end{liss}}

\newarrow{dashto}{}{dash}{}{dash}>

\newarrow{mapsto}|--->

\newarrow{dotsto}|...>

\begin{document}

\renewcommand{\a}{\alpha}
\renewcommand{\b}{\beta}
\newcommand{\g}{\gamma}
\renewcommand{\d}{\delta}
\newcommand{\e}{\epsilon}
\newcommand{\ve}{\varepsilon}
\newcommand{\z}{\zeta}
\renewcommand{\t}{\theta}
\renewcommand{\l}{\lambda}
\renewcommand{\k}{\varkappa}
\newcommand{\m}{\mu}
\newcommand{\n}{\nu}
\renewcommand{\r}{\rho}
\newcommand{\vr}{\varrho}
\newcommand{\s}{\sigma}
\newcommand{\vp}{\varphi}
\renewcommand{\o}{\omega}

\renewcommand{\Re}{\text{\rm Re }}

\newcommand{\G}{\Gamma}
\newcommand{\D}{\Delta}
\newcommand{\T}{\Theta}
\renewcommand{\L}{\Lambda}
\renewcommand{\P}{\Pi}
\newcommand{\Si}{\Sigma}
\renewcommand{\O}{\Omega}

\newcommand{\Up}{\Upsilon}

\renewcommand{\AA}{{\mathcal A}}
\newcommand{\BB}{{\mathcal B}}
\newcommand{\CC}{{\mathcal C}}
\newcommand{\DD}{{\mathcal D}}
\newcommand{\EE}{{\mathcal E}}
\newcommand{\FF}{{\mathcal F}}
\newcommand{\GG}{{\mathcal G}}
\newcommand{\HH}{{\mathcal H}}
\newcommand{\II}{{\mathcal I}}
\newcommand{\JJ}{{\mathcal J}}
\newcommand{\KK}{{\mathcal K}}
\newcommand{\LL}{{\mathcal L}}
\newcommand{\MM}{{\mathcal M}}
\newcommand{\NN}{{\mathcal N}}
\newcommand{\OO}{{\mathcal O}}
\newcommand{\PP}{{\mathcal P}}
\newcommand{\QQ}{{\mathcal Q}}
\newcommand{\RR}{{\mathcal R}}
\renewcommand{\SS}{{\mathcal S}}
\newcommand{\TT}{{\mathcal T}}
\newcommand{\UU}{{\mathcal U}}
\newcommand{\VV}{{\mathcal V}}
\newcommand{\WW}{{\mathcal W}}
\newcommand{\XX}{{\mathcal X}}
\newcommand{\YY}{{\mathcal Y}}
\newcommand{\ZZ}{{\mathcal Z}}

\renewcommand{\aa}{{\mathbb{A}}}
\newcommand{\bb}{{\mathbb{B}}}
\newcommand{\cc}{{\mathbb{C}}}
\newcommand{\dd}{{\mathbb{D}}}
\newcommand{\ee}{{\mathbb{E}}}
\newcommand{\ff}{{\mathbb{F}}}
\renewcommand{\gg}{{\mathbb{G}}}
\newcommand{\hh}{{\mathbb{H}}}
\newcommand{\ii}{{\mathbb{I}}}
\newcommand{\jj}{{\mathbb{J}}}
\newcommand{\kk}{{\mathbb{K}}}
\renewcommand{\ll}{{\mathbb{L}}}
\newcommand{\mm}{{\mathbb{M}}}
\newcommand{\nn}{{\mathbb{N}}}
\newcommand{\oo}{{\mathbb{O}}}
\newcommand{\pp}{{\mathbb{P}}}
\newcommand{\qq}{{\mathbb{Q}}}
\newcommand{\rr}{{\mathbb{R}}}
\renewcommand{\ss}{{\mathbb{S}}}
\newcommand{\ttt}{{\mathbb{T}}}
\newcommand{\uu}{{\mathbb{U}}}
\newcommand{\vv}{{\mathbb{V}}}
\newcommand{\ww}{{\mathbb{W}}}
\newcommand{\xx}{{\mathbb{X}}}
\newcommand{\yy}{{\mathbb{Y}}}
\newcommand{\zz}{{\mathbb{Z}}}

\newcommand{\AAA}{{\mathbf{A}}}
\newcommand{\BBB}{{\mathbf{B}} }
\newcommand{\CCC}{{\mathbf{C}} }
\newcommand{\DDD}{{\mathbf{D}} }
\newcommand{\EEE}{{\mathbf{E}} }
\newcommand{\FFF}{{\mathbf{F}} }
\newcommand{\GGG}{{\mathbf{G}}}
\newcommand{\HHH}{{\mathbf{H}}}
\newcommand{\III}{{\mathbf{I}}}
\newcommand{\JJJ}{{\mathbf{J}}}
\newcommand{\KKK}{{\mathbf{K}}}
\newcommand{\LLL}{{\mathbf{L}}}
\newcommand{\MMM}{{\mathbf{M}}}
\newcommand{\NNN}{{\mathbf{N}}}
\newcommand{\OOO}{{\mathbf{O}}}
\newcommand{\PPP}{{\mathbf{P}}}
\newcommand{\QQQ}{{\mathbf{Q}}}
\newcommand{\RRR}{{\mathbf{R}}}
\newcommand{\SSS}{{\mathbf{S}}}
\newcommand{\TTT}{{\mathbf{T}}}
\newcommand{\UUU}{{\mathbf{U}}}
\newcommand{\VVV}{{\mathbf{V}}}
\newcommand{\WWW}{{\mathbf{W}}}
\newcommand{\XXX}{{\mathbf{X}}}
\newcommand{\YYY}{{\mathbf{Y}}}
\newcommand{\ZZZ}{{\mathbf{Z}}}

\newcommand{\gA}{{\mathfrak{A}}}
\newcommand{\gB}{{\mathfrak{B}}}
\newcommand{\gC}{{\mathfrak{C}}}
\newcommand{\gD}{{\mathfrak{D}}}
\newcommand{\gE}{{\mathfrak{E}}}
\newcommand{\gF}{{\mathfrak{F}}}
\newcommand{\gG}{{\mathfrak{G}}}
\newcommand{\gH}{{\mathfrak{H}}}
\newcommand{\gI}{{\mathfrak{I}}}
\newcommand{\gJ}{{\mathfrak{J}}}
\newcommand{\gK}{{\mathfrak{K}}}
\newcommand{\gL}{{\mathfrak{L}}}
\newcommand{\gM}{{\mathfrak{M}}}
\newcommand{\gN}{{\mathfrak{N}}}
\newcommand{\gO}{{\mathfrak{O}}}
\newcommand{\gP}{{\mathfrak{P}}}
\newcommand{\gQ}{{\mathfrak{Q}}}
\newcommand{\gR}{{\mathfrak{R}}}
\newcommand{\gS}{{\mathfrak{S}}}
\newcommand{\gT}{{\mathfrak{T}}}
\newcommand{\gU}{{\mathfrak{U}}}
\newcommand{\gV}{{\mathfrak{V}}}
\newcommand{\gW}{{\mathfrak{W}}}
\newcommand{\gX}{{\mathfrak{X}}}
\newcommand{\gY}{{\mathfrak{Y}}}
\newcommand{\gZ}{{\mathfrak{Z}}}

\newcommand{\gota}{{\mathfrak{a}}}
\newcommand{\gotb}{{\mathfrak{b}}}
\newcommand{\gotc}{{\mathfrak{c}}}
\newcommand{\gotd}{{\mathfrak{d}}}
\newcommand{\gote}{{\mathfrak{e}}}
\newcommand{\gotf}{{\mathfrak{f}}}
\newcommand{\gotg}{{\mathfrak{g}}}
\newcommand{\goth}{{\mathfrak{h}}}
\newcommand{\goti}{{\mathfrak{i}}}
\newcommand{\gotj}{{\mathfrak{j}}}
\newcommand{\gotk}{{\mathfrak{k}}}
\newcommand{\gotl}{{\mathfrak{l}}}
\newcommand{\gotm}{{\mathfrak{m}}}
\newcommand{\gotn}{{\mathfrak{n}}}
\newcommand{\goto}{{\mathfrak{o}}}
\newcommand{\gotp}{{\mathfrak{p}}}
\newcommand{\gotq}{{\mathfrak{q}}}
\newcommand{\gotr}{{\mathfrak{r}}}
\newcommand{\gots}{{\mathfrak{s}}}
\newcommand{\gott}{{\mathfrak{t}}}
\newcommand{\gotu}{{\mathfrak{u}}}
\newcommand{\gotv}{{\mathfrak{v}}}
\newcommand{\gotw}{{\mathfrak{w}}}
\newcommand{\gotx}{{\mathfrak{x}}}
\newcommand{\goty}{{\mathfrak{y}}}
\newcommand{\gotz}{{\mathfrak{z}}}



\newcommand{\kkrest}{\begin{picture}(14,14)
\put(00,04){\line(1,0){14}}
\put(00,02){\line(1,0){14}}
\put(06,-4){\line(0,1){14}}
\put(08,-4){\line(0,1){14}}
\end{picture}     }

\newcommand{\krest}{~\kkrest~}


\newcommand{\grd}{{\text{\rm grd}}}
\newcommand{\id}{\text{id}}
\newcommand{\Tb}{\text{ \rm Tb}}
\newcommand{\Log}{\text{\rm Log }}
\newcommand{\Wh}{\text{\rm Wh }}
\newcommand{\Ker}{\text{\rm Ker }}
\newcommand{\Ext}{\text{\rm Ext}}
\newcommand{\Hom}{\text{\rm Hom}}
\newcommand{\diam}{\text{\rm diam}}
\newcommand{\Homb}{\text{\rm Hom}b}
\newcommand{\Lg}{\text{\rm Lg }}
\newcommand{\ind}{\text{\rm ind}}
\newcommand{\rk}{\text{\rm rk }}
\renewcommand{\Im}{\text{\rm Im }}
\newcommand{\supp}{\text{\rm supp }}
\newcommand{\Int}{\text{\rm Int }}
\newcommand{\grad}{\text{\rm grad}}
\newcommand{\Fix}{\text{\rm Fix}}
\newcommand{\Exp}{\text{\rm Exp}}
\newcommand{\Per}{\text{\rm Per}}
\newcommand{\TL}{\text{\rm TL}}
\newcommand{\Id}{\text{\rm Id}}
\newcommand{\Vect}{\text{\rm Vect}}
\newcommand{\vvol}{\text{\rm vol}}
\newcommand{\Mat}{\text\rm Mat}
\newcommand{\Tub}{\text{\rm Tub}}
\newcommand{\Imm}{\text{\rm Im}}
\newcommand{\tn}{\text{\rm t.n.}}
\newcommand{\card}{\text{\rm card }}
\newcommand{\GL}{\text{\rm GL }}

\newcommand{\track}{\text{\rm Track}}
\newcommand{\sgn}{\text{\rm sgn}}
\newcommand{\Arctg}{\text{\rm Arctg }}
\newcommand{\tg}{\text{\rm tg }}
\newcommand{\Arcsin}{\text{\rm Arcsin }}


\newcommand{\bere}{\begin{rema}}
\newcommand{\bede}{\begin{defi}}

\renewcommand{\beth}{\begin{theo}}
\newcommand{\bele}{\begin{lemm}}
\newcommand{\bepr}{\begin{prop}}
\newcommand{\beeq}{\begin{equation}}
\newcommand{\bega}{\begin{gather}}
\newcommand{\begaa}{\begin{gather*}}
\newcommand{\been}{\begin{enumerate}}

\newcommand{\bedee}{\begin{defii}}
\newcommand{\bethh}{\begin{theoo}}
\newcommand{\belee}{\begin{lemmm}}
\newcommand{\beprr}{\begin{propp}}

\newcommand{\beco}{\begin{coro}}

\newcommand{\beal}{\begin{aligned}}

\newcommand{\enre}{\end{rema}}

\newcommand{\enco}{\end{coro}}
\newcommand{\enpr}{\end{prop}}
\newcommand{\enth}{\end{theo}}
\newcommand{\enle}{\end{lemm}}
\newcommand{\enen}{\end{enumerate}}
\newcommand{\enga}{\end{gather}}
\newcommand{\engaa}{\end{gather*}}
\newcommand{\eneq}{\end{equation}}
\newcommand{\enal}{\end{aligned}}

\newcommand{\bq}{\begin{equation}}
\newcommand{\bqq}{\begin{equation*}}


\renewcommand{\leq}{\leqslant}
\renewcommand{\geq}{\geqslant}
\newcommand{\vphi}{\varphi}
\newcommand{\vide}{  \emptyset  }
\newcommand{\bu}{\bullet}
\newcommand{\pfff}{\pitchfork}
\newcommand{\mx}{\mbox}

\newcommand{\mxx}[1]{\quad\mbox{#1}\quad}
 \newcommand{\mxxx}[1]{\hspace{0.1cm}\mbox{#1} \quad  }
\newcommand{\wi}{\widetilde}

\newcommand{\ove}{\overline}
\newcommand{\unde}{\underline}
\newcommand{\ptf}{\pitchfork}

\newcommand{\emp}{\emptyset}
\newcommand{\wh}{\widehat}

\newcommand{\sub}{ Subsection}

\newcommand{\lc}{\lceil}
\newcommand{\rc}{\rceil}
\newcommand{\sps}{\supset}

\newcommand{\sm}{\setminus}
\newcommand{\ems}{\varnothing}
\newcommand{\sbs}{\subset}

\newcommand{\subs}{\subsection}
\newcommand{\ity}{\infty}


\newcommand{\GC}{\GG\CC}
\newcommand{\GCT}{\GG\CC\TT}
\newcommand{\GT}{\GG\TT}

\newcommand{\GA}{\GG\AA}
\newcommand{\GRP}{\GG\RR\PP}

\newcommand{\GgC}{\GG\gC}
\newcommand{\GgCC}{\GG\gC\CC}

\newcommand{\GgCT}{\GG\gC\TT}

\newcommand{\GgCY}{\GG\gC\YY}
\newcommand{\GgCYT}{\GG\gC\YY\TT}

\newcommand{\GCCT}{\GG\gC\CC\TT}
\newcommand{\GCC}{\GG\gC\CC}

\newcommand{\GKSC}{\GG\KK\SS\gC}
\newcommand{\GKS}{\GG\KK\SS}


\newcommand{\dstr}[1]
{
\TT_{#1}
}

\newcommand{\strr}[3]
{{#1}^{\displaystyle\twoheadrightarrow}_{[{#2},{#3}]}}

\newcommand{\str}[1]{{#1}^{\displaystyle\twoheadrightarrow}}

\newcommand{\stind}[3]
{{#1}^{\displaystyle\rightsquigarrow}_{[{#2},{#3}]}}

\newcommand{\st}[1]{\overset{\rightsquigarrow}{#1}}
\newcommand{\bst}[1]{\overset{\displaystyle\rightsquigarrow}
\to{\boldkey{#1}}}

\newcommand{\stexp}[1]{{#1}^{\rightsquigarrow}}
\newcommand{\bstexp}[1]{{#1}^{\displaystyle\rightsquigarrow}}

\newcommand{\bstind}[3]{{\boldkey{#1}}^{\displaystyle\rightsquigarrow}_
{[{#2},{#3}]}}
\newcommand{\bminstind}[3]{\stind{({\boldkey{-}\boldkey{#1}})}{#2}{#3}}

\newcommand{\ST}{\stexp}

\newcommand{\stv}{\stexp {(-v)}}
\newcommand{\stu}{\stexp {(-u)}}
\newcommand{\stw}{\stexp {(-w)}}

\newcommand{\strv}[2]{\stind {(-v)}{#1}{#2}}
\newcommand{\strw}[2]{\stind {(-w)}{#1}{#2}}
\newcommand{\stru}[2]{\stind {(-u)}{#1}{#2}}

\newcommand{\stvv}[2]{\stind {v}{#1}{#2}}
\newcommand{\stuu}[2]{\stind {u}{#1}{#2}}
\newcommand{\stww}[2]{\stind {w}{#1}{#2}}

\newcommand{\vovo}{\stexp {(-v0)}}
\newcommand{\vov}{\stexp {(-v1)}}

\newcommand{\fl}[1]{{#1}\!\da}
\newcommand{\fll}[1]{({#1}\!\da)}

\newcommand{\vflesh}{\fl{v}}
\newcommand{\wflesh}{\fl{w}}

\newcommand{\vfllesh}{\fll{v}}
\newcommand{\wfllesh}{\fll{w}}

\newcommand{\vvfl}{\wi v\!\da}
\newcommand{\wivflesh}{\wi v\!\da}

\newcommand{\RA}{\Rightarrow}
\newcommand{\LA}{\Leftarrow}
\newcommand{\RLA}{\Leftrightarrow}

\newcommand{\LRA}{\Leftrightarrow}

\newcommand{\lau}[1]{{\xleftarrow{#1}}}

\newcommand{\rau}[1]{{\xrightarrow{#1}}}
\newcommand{\rad}[1]{ {\xrightarrow[#1]{}} }

\newcommand{\da}{\downarrow}


\newcommand{\vecm}{\Vect^1_0(M)}
\newcommand{\vecbn}{\Vect^1_b(N)}
\newcommand{\vecw}{\Vect^1(W)}
\newcommand{\vecrm}{\Vect^1_b(\RRR^m)}
\newcommand{\vecbm}{\Vect^1_b(M)}
\newcommand{\ver}{\text{\rm Vect}^1(\RRR^ n)}
\newcommand{\verr}{\text{\rm Vect}^1_0(\RRR^ n)}
\newcommand{\hrrr}{\text{\rm Vect}^1(M)}
\newcommand{\vemm}{\text{\rm Vect}^1_0(M)}

\newcommand{\vem}{\text{{\rm Vectt}}(M)}
\newcommand{\vebK}{\text{{\rm Vectt}}(B,K)}
\newcommand{\vemK}{\text{{\rm Vectt}}(M,K)}
\newcommand{\vemc}{\text{{\rm Vectt}}_c(M)}
\newcommand{\vemQ}{\text{{\rm Vectt}}(M,Q)}

\newcommand{\vectt}[1]{\text{{\rm Vectt}}(#1)}

\newcommand{\vecsmo}[1]{{\text{\rm Vect}}^\infty (#1)}

\newcommand{\vew}{\text{\rm Vect}^1 (W,\bot)}

\newcommand{\downnorm}{\text{\rm Vect}^1_\bot (W)}

\newcommand{\upnorm}{\text{\rm Vect}^1_\top (W)}

\newcommand{\normm }{\text{\rm Vect}^1_N (W)}

\newcommand{\veww}{\text{\rm Vect}^1 (W)}


\newcommand{\tens}[1]{\underset{#1}{\otimes}}

\newcommand{\starr}[1]{\underset{#1}{*}}

\newcommand{\Lxxxi}{\wh L_{\bar\xi}}

\newcommand{\Lxxi}{\L_{ [\xi] }}

\newcommand{\LLLxi}{\wh\LL_\xi}

\newcommand{\Lbarxi}{\wh \L_{\bar\xi}}

\newcommand{\RRxi}{\wh \RR_\xi}
\newcommand{\LLxi}{\wi \L_\xi}

\newcommand{\Lx}{\L_{(\xi)}}
\newcommand{\Lxi}{{\wh \L}_\xi}
\newcommand{\Leta}{{\wh \L}_\eta}

\newcommand{\lL}{\wh{\wh L}}

\newcommand{\Rxi}{{\ove R}_\xi}
\newcommand{\Nxi}{{\ove N}_\xi}
\newcommand{\Rcxi}{{\bar R}_\xi^c}

\newcommand{\sil}{ S^{-1}\L }
\newcommand{\kil}{\ove{K}_1(\L)}
\newcommand{\killl}{\ove{K}_1(\wh\L)}
\newcommand{\kisl}{\ove{K}_1(S^{-1}\L )}

\newcommand{\klxi}{K_1(\Lxi)}
\newcommand{\kklxi}{\ove{K_1}(\Lxi)}

\newcommand{\popo}{\tens{\L}\Lxi}
\newcommand{\popom}{\tens{\L^-}\Lxi^-}


\newcommand{\amk}{\AA^{(m)}_k}
\newcommand{\amkm}{\AA^{(m)}_{k-1}}
\newcommand{\bmk}{\BB^{(m)}_k}
\newcommand{\bmkm}{\BB^{(m)}_k}

\newcommand{\tivm}{\wi V^-}

\newcommand{\vk}{V_{\langle k\rangle}^-}
\newcommand{\tivkm}{\wi V_{\langle k-1\rangle}^-}
\newcommand{\tivk}{\wi V_{\langle k\rangle}^-}

\newcommand{\vkm}{V_{\langle k-1\rangle}^-}
\newcommand{\vkp}{V_{\langle k+1\rangle}^-}

\newcommand{\hvk}{\wh V_{\langle k\rangle}^-}

\newcommand{\hvkm}{\wh V_{\langle k-1\rangle}^-}
\newcommand{\hvkp}{\wh V_{\langle k+1\rangle}^-}

\newcommand{\vkvk}{V_{\prec k\succ}^-}

\newcommand{\vkvkm}{V_{\prec k-1\\succ}^-}

\newcommand{\tivkvk}{\wi V_{\prec k\succ}^-}

\newcommand{\tivkvkm}{\wi V_{\prec k-1 \succ}^-}

\newcommand{\vvvbs}{V_b^{( s)    }}
\newcommand{\vvvas}{V_a^{( s)    }}
\newcommand{\vvvbsm}{V_b^{( s-1)    }}
\newcommand{\vvvasm}{V_a^{( s-1)    }}

\newcommand{\vvbsm}{V_b^{( s-1)    }}
\newcommand{\vvasm}{V_a^{( s-1)    }}
\newcommand{\vvbs}{V_b^{[s]}    }
\newcommand{\vvas}{V_a^{[ s]    }}

\newcommand{\factor}{\vvbs / \vvbsm}
\newcommand{\factora}{\vvas / \vvasm}

\newcommand{\vvksm}{\wi V_k^{( s-1)   }}
\newcommand{\vvks}{\wi V_k^{[s]}    }

\newcommand{\vvkmsm}{\wi V_{k-1}^{( s-1)   }}
\newcommand{\vvkms}{\wi V_{k-1}^{[s]}    }

\newcommand{\fac}{\vvks / \vvksm}
\newcommand{\facm}{\vvkms / \vvkmsm}

\newcommand{\vbsm}{V_b^{\{\leq s-1\}}    }
\newcommand{\vasm}{V_a^{\{\leq s-1\}}    }
\newcommand{\vbs}{V_b^{\{\leq s\}}    }
\newcommand{\vas}{V_a^{\{\leq s\}}    }

\newcommand{\wivksm}{\wi V_k^{\{\leq s-1\}}    }
\newcommand{\wivkmsm}{\wi V_{k-1}^{\{\leq s-1\}}    }
\newcommand{\wivks}{\wi V_k^{\{\leq s\}}    }
\newcommand{\wivkms}{\wi V_{k-1}^{\{\leq s\}}    }

\newcommand{\Vbsm}{V_b^{[\leq s-1]}(\d)    }
\newcommand{\Vasm}{V_a^{[\leq s-1]}(\d)    }
\newcommand{\Vbs}{V_b^{[\leq s]}(\d)    }
\newcommand{\Vas}{V_a^{[\leq s]}(\d)    }

\newcommand{\vass}{V_{a_{s+1}}}

\newcommand{\vbkm}{V_b^{\{\leq k-1\}}    }
\newcommand{\vakm}{V_a^{\{\leq k-1\}}    }
\newcommand{\vbk}{V_b^{\{\leq k\}}    }
\newcommand{\vak}{V_a^{\{\leq k\}}    }

\newcommand{\Vbkm}{V_b^{[\leq k-1]}(\d)    }
\newcommand{\Vakm}{V_a^{[\leq k-1]}(\d)    }
\newcommand{\Vbk}{V_b^{[\leq k]}(\d)    }
\newcommand{\Vak}{V_a^{[\leq s]}(\d)    }


\newcommand{\dow}{\pr_0 W}

\newcommand{\daw}{\pr_1 W}

\newcommand{\hdaw}{\wh{\pr_1 W}}

\newcommand{\tipwk}{(\pr_1 \wi W)^{\{\leq k\}}}

\newcommand{\timwk}{(\pr_0 \wi W)^{\{\leq k\}}}

\newcommand{\tipwkm}{(\pr_1 \wi W)^{\{\leq k-1\}}}

\newcommand{\pws}{(\pr_1 W)^{\{\leq s\}}}

\newcommand{\hpws}{(\wh{\pr_1 W})^{\{\leq s\}}}

\newcommand{\hpwsm}{(\wh{\pr_1 W})^{\{\leq {s-1}\}}}

\newcommand{\mws}{(\pr_0 W)^{\{\leq s\}}}

\newcommand{\pwk}{(\pr_1 W)^{\{\leq k\}}}

\newcommand{\mwk}{(\pr_0 W)^{\{\leq k\}}}

\newcommand{\pwkm}{(\pr_1 W)^{\{\leq k-1\}}}

\newcommand{\pwsm}{(\pr_1 W)^{\{\leq s-1\}}}

\newcommand{\pwkmm}{(\pr_1 W)^{\{\leq k-2\}}}

\newcommand{\mwkmm}{(\pr_0 W)^{\{\leq k-2\}}}

\newcommand{\mwkm}{(\pr_0 W)^{\{\leq k-1\}}}

\newcommand{\mwsm}{(\pr_0 W)^{\{\leq s-1\}}}

\newcommand{\mwkp}{(\pr_0 W)^{\{\leq k+1\}}}

\newcommand{\dwmok}{\daw^{( k)}}

\newcommand{\dwmokp}{\daw^{(k+1)}}

\newcommand{\dwmokm}{\daw^{(k-1)}}
\newcommand{\dwmokmm}{\daw^{(k-2)}}


\newcommand{\moi}[1]{\MM^{(0)}_{#1}}
\newcommand{\moii}[1]{\MM^{(1)}_{#1}}

\newcommand{\Wal}{W_{[a,\l]}}

\newcommand{\Wlm}{W_{[\l,\m]}}
\newcommand{\Wlb}{W_{[\l,b]}}

\newcommand{\Wam}{W_{[a,\m]}}

\newcommand{\Wall}{W_{[a,\l']}}

\newcommand{\Wkr}{W^{\circ}}
\newcommand{\wkr}{W^{\circ}}

\newcommand{\wa}[2]{ W_{[a_{#1}, a_{#2}]}}

\newcommand{\waa}[1]{ W_{[a, a_{#1}]}}

\newcommand{\Wa}[2]{ W_{[{#1}, {#2}]}}

\newcommand{\WS}[1]{ W^{\{\leq {#1}\}}}

\newcommand{\ws}{\WS {s}}

\newcommand{\wsm}{\WS {s-1}}

\newcommand{\wsmm}{\WS {s-2}}

\newcommand{\wk}{\WS {k}}

\newcommand{\wkm}{\WS {k-1}}

\newcommand{\wkmm}{\WS {k-2}}

\newcommand{\wsn}{ W^{[\leq s]}(\nu)}

\newcommand{\wsmn}{ W^{[\leq s-1]}(\nu)}

\newcommand{\wsk}{ W^{[\leq k]}(\nu)}

\newcommand{\Wmok}{W^{\prec k\succ}}
\newcommand{\wmok}{W^{\langle k\rangle}}

\newcommand{\wmokp}{W^{\langle k+1\rangle}}

\newcommand{\wmokm}{W^{\langle k-1\rangle}}
\newcommand{\wmokmm}{W^{\langle k-2\rangle}}

\newcommand{\wmos}{W^{\langle s\rangle}}

\newcommand{\wmosm}{W^{\langle s-1\rangle}}

\newcommand{\wmoo}{W^{\langle 0\rangle}}

\newcommand{\wwk}{\( \wmok , \wmokm \)}

\newcommand{\wwkp}{\( \wmokp , \wmok \)}

\newcommand{\wwkm}{\( \wmokm , \wmokmm \)}

\newcommand{\wws}{\bigg( \wmos , \wmosm \bigg)}

\newcommand{\wasn}{W^{\lc s\rc}(\nu)}

\newcommand{\wakn}{W^{\lc k\rc}(\nu)}

\newcommand{\hwm}{H_*\( \wmok, \wmokm\)}

\newcommand{\hkwm}{H_k\( \wmok, \wmokm\)}

\newcommand{\hwmp}{H_*\( \wmokp, \wmok\)}

\newcommand{\hkwmm}{H_k\( \wmokm, \wmokmm\)}

\newcommand{\hkwmp}{H_k\( \wmokp, \wmok\)}

\newcommand{\tiws}{\wi W^{\{\leq s\}}}

\newcommand{\tiwk}{\wi W^{\{\leq k\}}}

\newcommand{\tiwsm}{\wi W^{\{\leq s-1\}}}

\newcommand{\tiwkm}{\wi W^{\{\leq k-1\}}}

\newcommand{\wkwk}{ W^{[ k]}}
\newcommand{\tiwmok}{\wi W^{\prec k\succ}}
\newcommand{\womok}{W_0^{\prec k\succ}}

\newcommand{\Womok}{W_0^{\langle k\rangle}}

\newcommand{\Wmokp}{W^{\langle k+1\rangle}}

\newcommand{\Wmokm}{W^{\langle k-1\rangle}}
\newcommand{\Wmokmm}{W^{\langle k-2\rangle}}

\newcommand{\Wmos}{W^{\langle s\rangle}}

\newcommand{\Wmosm}{W^{\langle s-1\rangle}}

\newcommand{\Wmoo}{W^{\langle 0\rangle}}

\newcommand{\hWmok}{\wh W^{\langle k\rangle}}
\newcommand{\hWmokm}{\wh W^{\langle k-1\rangle}}

\newcommand{\tiWmok}{\wi W^{\langle k\rangle}}
\newcommand{\tiWmokm}{\wi W^{\langle k-1\rangle}}



\newcommand{\talp}{twisted Alexander polynomial}
\newcommand{\tnh}{twisted Novikov homology}

\newcommand{\ifff}{if and only if}

\newcommand{\orial}{oriented almost transverse}
\renewcommand{\th}{therefore}
\newcommand{\at}{almost~ transverse}
\newcommand{\ata}{almost~ transversality~ condition}
\newcommand{\gr}{gradient}
\newcommand{\Mf}{Morse function}
\newcommand{\iis}{it is sufficient}
\newcommand{\sut}{~such~that~}
\newcommand{\sufsm}{~sufficiently~ small}
\newcommand{\sufla}{~sufficiently~ large}
\newcommand{\sufcl}{~sufficiently~ close}
\newcommand{\wrt}{~with respect to}
\newcommand{\ho}{homomorphism}
\newcommand{\iso}{isomorphism}
\newcommand{\rgr}{Riemannian gradient}
\newcommand{\riemm}{Riemannian metric}

\newcommand{\trasp}{trajectory starting at a point of}
\newcommand{\trasps}{trajectories starting at a point of}

\newcommand{\ma}{manifold}
\newcommand{\nei}{neighbourhood}
\newcommand{\dfm}{diffeomorphism}

\newcommand{\vf}{vector field}

\newcommand{\vfs}{vector fields}

\newcommand{\fe}{for every}

\newcommand{\tr}{~trajectory }

\newcommand{\grs}{~gradients}
\newcommand{\trs}{~trajectories}

\newcommand{ \co}{~cobordism}
\newcommand{
\sma}{submanifold}
\newcommand{
\hos}{~homomorphisms}
\newcommand{
\Th}{~Therefore}

\newcommand{
\tthen}{\text \rm ~then}

\newcommand{
\wwe}{\text \rm ~we  }
\newcommand{
\hhave}{\text \rm ~have}
\newcommand{
\eevery}{\text \rm ~every}

\newcommand{\noconf}{~no~confusion~is~possible}

\newcommand{\ATA}{almost~ transversality~ tondition}
\newcommand{\cob}{~cobordism}

\newcommand{\hot}{homotopy}

\newcommand{\emem}{elementary modification}
\newcommand{\emems}{elementary modifications}

\newcommand{\TA}{transversality condition}

\newcommand{\hog}{homology}

\newcommand{\cog}{cohomology}

\newcommand{\wat}{ We shall assume that}

\newcommand{\sclv}{sufficiently close to $v$ in $C^0$-topology}

\newcommand{\cf}{continuous function }

\newcommand{\heq}{homotopy equivalence}

\newcommand{\heeq}{homology equivalence}

\newcommand{\eg}{exponential growth}

\newcommand{\nics}{Novikov incidence coefficients}
\newcommand{\nic}{Novikov incidence coefficient}

\newcommand{\negc}{Novikov exponential growth conjecture}

\newcommand{\mc}{Morse Complex   }

\newcommand{\mas}{manifolds   }

\newcommand{\nc}{Novikov Complex   }

\newcommand{\glvf}{gradient-like vector field}

\newcommand{\glvfs}{gradient-like vector fields}

\newcommand{\fg}{finitely generated   }

\newcommand{\she}{simple~homotopy~equivalence}

\newcommand{\sht}{simple~homotopy~type}

\newcommand{\ta}{transversality condition}

\newcommand{\cpc}{convex polyhedral cone}
\newcommand{\rcpc}{rational convex polyhedral cone}

\newcommand{\mnp}{Morse-Novikov pair}

\newcommand{\rp}{rationality property}
\newcommand{\wvf}{Whitney vector field}

\newcommand{\egp}{exponential growth property}

\newcommand{\lzf}{Lefschetz zeta function}

\newcommand{\babs}{by abuse of notation}
\newcommand{\su}{subsection}
\newcommand{\Prop}{\text{Proposition}}

\newcommand{\aand}{\quad\text{and}\quad}
\newcommand{\wwhere}{\quad\text{where}\quad}
\newcommand{\ffor}{\quad\text{for}\quad}
\newcommand{\iif}{\quad\text{if}\quad}
\newcommand{\iiif}{~\text{if}~}

\newcommand{\eqi}{equivalence}

\newcommand{\mfcobv}{~ Let $\fcob$
be a Morse function on
 a cobordism $W$ and $v$
 be an $f$-gradient. ~}

\newcommand{\mfcob}{~ Let $\fcob$
be a Morse function on
 a cobordism $W$}

\newcommand{\mfcobvat}{~ Let $\fcob$
be a Morse function on
 a cobordism $W$ and $v$
  an almost transverse $f$-gradient. ~}

\newcommand{\msf}{Morse-Smale filtration}

\newcommand{\fbfg}{ free based finitely generated }

\newcommand{\tap}
{twisted Alexander invariant}


\newcommand{\vaa}{\mathscr A_k}

\newcommand{\qaa}{\mathscr Q_k}

\newcommand{\tret}{{\frac 13}}
\newcommand{\dvet}{{\frac 23}}
\newcommand{\polt}{{\frac 32}}
\newcommand{\polo}{{\frac 12}}

\newcommand{\bv}{B(v,\d)}

\newcommand{\ti}{\times}

\newcommand{\FR}{{\mathcal{F}}r}
\newcommand{\gt}{{\mathcal{G}}t}

\newcommand{\en}{enumerate}

\newcommand{\Prf}{{\it Proof.\quad}}
\newcommand{\prf}{{\it Proof:\quad}}

\newcommand{\nr}{\Vert}
\newcommand{\smo}{C^{\infty}}

\newcommand{\fpr}[2]{{#1}^{-1}({#2})}
\newcommand{\sdvg}[3]{\widehat{#1}_{[{#2},{#3}]}}
\newcommand{\disc}[3]{B^{({#1})}_{#2}({#3})}
\newcommand{\Disc}[3]{D^{({#1})}_{#2}({#3})}
\newcommand{\desc}[3]{B_{#1}(\leq{#2},{#3})}
\newcommand{\Desc}[3]{D_{#1}(\leq{#2},{#3})}
\newcommand{\komp}[3]{{\bold K}({#1})^{({#2})}({#3})}
\newcommand{\Komp}[3]{\big({\bold K}({#1})\big)^{({#2})}({#3})}

\newcommand{\ran}{\{(A_\lambda , B_\lambda)\}_{\lambda\in\Lambda}}
\newcommand{\rran}{\{(A_\lambda^{(s)},
 B_\lambda^{(s)}  )\}_{\lambda\in\Lambda, 0\leq s\leq n }}
\newcommand{\brs}{\rran}
\newcommand{\rans}{\{(A_\sigma , B_\sigma)\}_{\sigma\in\Sigma}}

\newcommand{\fmin}{F^{-1}}
\newcommand{\vh}{\widehat{(-v)}}

\newcommand{\chart}{\Phi_p:U_p\to B^n(0,r_p)}
\newcommand{\atlas}{\{\Phi_p:U_p\to B^n(0,r_p)\}_{p\in S(f)}}
\newcommand{\flow}{{\VV}=(f,v, \UU)}

\newcommand{\Rn}{\bold R^n}
\newcommand{\Rk}{\bold R^k}

\newcommand{\fcob}{f:W\to[a,b]}

\newcommand{\phicob}{\phi:W\to[a,b]}

\newcommand{\crr}{p\in S(f)}
\newcommand{\nrv}{\Vert v \Vert}
\newcommand{\nrw}{\Vert w \Vert}
\newcommand{\nru}{\Vert u \Vert}

\newcommand{\obb}{\cup_{p\in S(f)} U_p}
\newcommand{\proob}{\Phi_p^{-1}(B^n(0,}

\newcommand{\indl}[1]{{\scriptstyle{\text{\rm ind}\leqslant {#1}~}}}
\newcommand{\inde}[1]{{\scriptstyle{\text{\rm ind}      =   {#1}~}}}
\newcommand{\indg}[1]{{\scriptstyle{\text{\rm ind}\geqslant {#1}~}}}

\newcommand{\obbi}{\cup_{p\in S_i(f)}}

\newcommand{\pr}{\partial}
\newcommand{\prx}[2]{\frac {\pr {#1}}{\pr x} ({#2})}
\newcommand{\pry}[2]{\frac {\pr {#1}}{\pr y} ({#2})}
\newcommand{\prz}[2]{\frac {\pr {#1}}{\pr z} ({#2})}
\newcommand{\przbar}[2]{\frac {\pr {#1}}{\pr \bar z} ({#2})}
\newcommand{\chape}[2]{\frac  {\pr {#1}}{\pr {#2}} }
\newcommand{\chapee}[2]{\frac  {\pr^2 {#1}}{\pr {#2}^2} }

\newcommand{\xit}{\tilde\xi_t}

\newcommand{\VODIN}{V_{1/3}}
\newcommand{\VDVA}{V_{2/3}}
\newcommand{\VM}{V_{1/2}}
\newcommand{\ddd}{\cup_{p\in S_i(F_1)} D_p(u)}
\newcommand{\dddmin}{\cup_{p\in S_i(F_1)} D_p(-u)}
\newcommand{\where}{\quad\text{\rm where}\quad}

\newcommand{\kr}[1]{{#1}^{\circ}}

\newcommand{\mods}{\vert s(t)\vert}
\newcommand{\exd}{e^{2(D+\alpha)t}}
\newcommand{\exmin}{e^{-2(D+\alpha)t}}

\newcommand{\intt}{[-\theta,\theta]}

\newcommand{\ffmin}{f^{-1}}

\newcommand{\vxi}{v\langle\vec\xi\rangle}

\newcommand{\qt}{\hfill\triangle}
\newcommand{\qs}{\hfill\square}

\newcommand{\pa}{\vskip0.1in}

\renewcommand{\(}{\big(}
\renewcommand{\)}{\big)}

\newcommand{\Vm}{V_\m}

\newcommand{\Vl}{V_\l}

\newcommand{\lccc}{\wh\L_{C}}

\newcommand{\ld}{\wh\L_{D}}

\newcommand{\udp}{{\displaystyle {\vartriangle}}}
\newcommand{\ddp}{{\displaystyle {\triangledown}}}

\newcommand{\Vv}{{\boldsymbol{v}}}

\newcommand{\hV}{\wh V}
\newcommand{\hHH}{\wh \HH}

\newcommand{\gama}[2]{\g({#1}, \tau_a({#2},{#1}); w )}

\newcommand{\gam}[2]{\g({#1}, \tau_0({#2},{#1}); w )}
\newcommand{\ga}[2]{\g({#1}, \tau({#2},{#1}); w )}

\newcommand{\mi}[3]{{#1}^{-1}\([{#2},{#3}]\)}

\newcommand{\fii}[2]{\mi {\phi}{a_{#1}}{a_{#2}} }

\newcommand{\fifi}[2]{\mi {\phi}{#1}{#2} }

\newcommand{\pf}[2]{\mi {\phi_1}{\a_{#1}}{\a_{#2}} }

\newcommand{\mf}[2]{\mi {\phi_0}{\b_{#1}}{\b_{#2}}}

\newcommand{\dqr}{\pr_- Q_r}

\newcommand{\ds}{\pr_s}

\newcommand{\dsm}{\pr_{s-1}}

\newcommand{\yz}{Y_k(v)\cup Z_k(v)}

\newcommand{\Gama}{{\nazad{ \Gamma}}}
\newcommand{\ug}[1]{\llcorner {#1} \lrcorner}
\newcommand{\npqv}{n(\bar p, \bar q; v)}
\newcommand{\fms}{f:M\to S^1   }

\newcommand{\nkpqv}{n_k(\bar p, \bar q; v)}

\newcommand{\GLT}{\GG lt}

\newcommand{\Trln}{{\text Trln}}

\newcommand{\Trlln}{{\text TrLn}}

\newcommand{\Tr}{{\text{\rm  Tr}}}
\newcommand{\TrL}{{\text TrL}}
\newcommand{\limdir}{\underset {\to}{\lim}}
\newcommand{\liminv}{\underset {\leftarrow}{\lim}}

\newcommand{\kom}[2]{ {#1}{#2}{ {#1}^{-1}} {{#2}^{-1}} }

\newcommand{\komm}[2]{ {#1}{#2}{ ({#1})^{-1}} {({#2})^{-1}} }
\newcommand{\kommm}[2]{ {#1}'{#2}'{ ({#1}'')^{-1}} {({#2}'')^{-1}} }

\newcommand{\Trll}{\TL'}

\newcommand{\cmd}{ C_*^\D( \wi M)}
\newcommand{\cmxi}{\wh C_*^\D( \wi M, \xi)}
\newcommand{\whgxi}{\wh {{\rm Wh}} (G,\xi)}

\newcommand{\vwdwp}{Vect(W,\pr_0W;P)}
\newcommand{\ewdwp}{\EE(W,\pr_0W)}
\newcommand{\ewdwo}{\EE(W_1,\pr_0W_1)}
\newcommand{\ewdwd}{\EE(W_2,\pr_0W_2)}

\newcommand{\kpr}{K_r^+}

\newcommand{\kmr}{K_r^-}

\newcommand{\kpd}{K_r^+(\d)}

\newcommand{\kmd}{K_r^-(\d)}

\newcommand{\addc}{\addtocontents{toc}{\protect\vspace{10pt}}}

\newcommand{\mxi}{M_\xi   }

\newcommand{\cmm}{C_*^\D(\wi M)}

\newcommand{\cvm}{C_*^\D(\wi V^-)}

\newcommand{\ey}{\wi E_*}
\newcommand{\eey}{\wi \EE_*}

\newcommand{\eky}{\wi E(k)_*}

\newcommand{\eti}{\wi{\wi\EE_*}}

\newcommand{\etik}{\wi{\wi\EE}_k}

\newcommand{\etikp}{\wi{\wi\EE}_{k+1}}

\newcommand{\ctiu}{\wi{\wi C}_*(u_1)}

\newcommand{\ctiuk}{\wi{\wi C}_k(u_1)}

\newcommand{\ctiv}{\wi{\wi C}_*(v)}

\newcommand{\ctiukm}{\wi{\wi C}_{k-1}(u_1)}

\newcommand{\scc}[1]{|{\scriptscriptstyle{#1}}}\newcommand{\rrr}{\{\wi r\}}

\newcommand{\tidow}{\pr_0 \wi W}
\newcommand{\tidaw}{\pr_1 \wi W}

\newcommand{\tivkp}{\wi V_{\langle k+1\rangle}^-}

\newcommand{\ur}[1]{\overset{\smallfrown}{#1}}

\newcommand{\dr}[1]{\underset{\smallsmile}{#1}}

\newcommand{\uUu}{\overset{\twoheadrightarrow}{u}}
\newcommand{\vVv}{\overset{\twoheadrightarrow}{v}}
\newcommand{\wWw}{\overset{\twoheadrightarrow}{w}}

\newcommand{\bfun}{{\bf 1}}

\newcommand{\ppmm}{{\scriptstyle{ \pm}}}

\newcommand{\dpm}{{\scriptstyle{ \pm}}}


\newcommand{\Lxim}{{\wh \L}^-_\xi}
\newcommand{\Lgxi}{{\wh \gL}_\xi}
\newcommand{\lxi}{{\bar \L}_\xi}

\newcommand{\Lc}{{\wh \L}_C}
\newcommand{\Lcm}{{\wh \L}_C^-}

\newcommand{\lLL}{\wh{\wh \L}}

\newcommand{\Xc}{{\wh X}_C}
\newcommand{\Xfaa}{{\wh X}_{(F_i,\vec\a)}}

\newcommand{\bs}{\boldsymbol}


\newcommand{\bikl}{\text{\rM (}}
\newcommand{\bikr}{\text{\rM )}}

\newcommand{\ck}{C_*^{(k)}}
\newcommand{\dk}{D_*^{(k)}}
\newcommand{\ckp}{C_*^{(k+1)}}
\newcommand{\dkp}{D_*^{(k+1)}}

\newcommand{\fmk}{F_m^{(k)}}
\newcommand{\fmkk}[1]{F_m^{(k)} ({#1}) }




\newcommand{\arnoldXC}{ V.I. Arnold,
\emph{Dynamics of intersections},
        Proceedings of a Conference
 in Honour of
J.Moser, edited by
 P.Rabinowitz and R.Hill,
 Academic Press,
 1990
pp. 77--84.  }

\newcommand{\arnoldXCprim}{ V.I. Arnold,
\emph{Dynamics of Complexity of Intersections},
Boletim SOc. Brasil. Mat. (N.S.),
1990, {\bff 21} (1), 1-10.}

\newcommand{\arnoldXCIII}{ V.I. Arnold,
\emph{Bounds for Milnor numbers of intersections
in holomorphic dynamical systems}, In:
Topological Methods in Modern Mathematics,
 Publish or Perish,
 1993,
pp. 379--390.}

\newcommand{\arnoldXCIV}{ V.I. Arnold,
\emph{Sur quelques probl\`emes de
 la th\'eorie des syst\`emes
dynamiques},
  Journal of the Julius Schauder center, {\bff 4}
 1994
pp. 209--225.  }

\newcommand{\armaz}
{M.Artin, B.Mazur,
\emph{On periodic points},
Annals of Math. {\bff 102} (1965), 82--99.
}

\newcommand{\baladi}
{V.Baladi,
\emph{Periodic orbits and dynamical spectra},
Ergodic theory and dynamical spectra,
{\bff 18}
(1998),
    255 - 292.
}

\newcommand{\chapman}
{ T.A.Chapman,
\emph{ Topological invariance of Whitehead torsion},
American J. of Math.
{\bff 96}
(1974),
    488 - 497
}

\newcommand{\CockroftSwan}[1]
{W.H.Cockroft, R.G.Swan, {#1}}

\newcommand{\twodim}
{
\emph{ On the homotopy type of 
2-dimensional complexes}
J.London Math. Soc, 
{\bff 11}, 1 (1961) 306-311. }

\newcommand{\conway}
{ J.H.Conway,
\emph{
 An enumeration of 
knots and links, and some 
of their algebraic properties},
Computational 
Problems in Abstract Algebra,
Pergamon Press, 
New York, 1970, pp. 329--358.}

\newcommand{\bhs}
{H.Bass, A.Heller, R.G.Swan,
\emph{The Whitehead group of a polynomial extension},
Inst. Hautes Etudes Sci. Publ. Math. {\bff 22} (1964),
61--79
}

\newcommand{\brlev}
{ W.Browder, J.Levine,
\emph{Fibering manifolds over the circle},
Comment. Math. Helv. {\bff 40} (1966),
153--160
}

\newcommand{\farran}{ M. Farber and A. A. Ranicki,
\emph{ The Morse-Novikov theory of circle-valued functions
and noncommutative localization,}

  Proc. 1998 Moscow Conference
for the 60th Birthday of S. P. Novikov, tr.
 Mat. Inst. Steklova, {\bf 225}, 1999,
381 -- 388.

\quad E-print:

dg-ga/9812122.
}

\newcommand{\farrell}
{F.T.Farrell,
\emph{The obstruction to fibering a manifold over a circle},
Indiana Univ.~J. {\bff 21} (1971), 315--346.
}

\newcommand{\farhsi}
{F.T.Farrell, W.-C.Hsiang,
\emph{A formula for $K_1R_\alpha[T]$},
Proc. Symp. Pure Math., Vol. {\bff 17} (1968), 192--218}

\newcommand{\fel}{  A.Fel'shtyn,
\emph{
Dynamical zeta functions,
Nielsen Theory and Reidemeister torsion,}
preprint ESI 539 of The Erwin Schr\"odinger International Institute for
Mathematical Physics,
(to appear in Memoirs of AMS)
(1998)
}

\newcommand{\franks}
{   J.Franks
\emph{Homology and dynamical systems},
CBMS Reg. Conf. vol. 49, AMS, Providence 1982.
}

\newcommand{\fuller}{  F.B.Fuller,
\emph{An index of fixed point type for periodic orbits},
Amer. J.Math
{\bff 89},
(1967)
133--148
}

\newcommand{\fried}{  D.Fried,
\emph{Homological Identities for closed orbits},
Inv. Math. {\bff 71}, (1983) 419--442.
}

\newcommand{\friedtwi}{  D.Fried,
\emph{Periodic points and twisted coefficients},
Lect. Notes in Math.,
{\bff 1007},
(1983)
261--293.
}

\newcommand{\friednewzeta}{  D.Fried,
\emph{Flow equivalence, hyperbolic systems and a new zeta function
for flows},
Comm. Math. Helv.,
{\bff 57},
(1982)
237--259.
}

\def\gnXCIV
{  R.Geoghegan, A.Nicas,
\emph{Trace and torsion in the theory of flows},
Topology,
{\bff 33},
(1994)
683--719
}

\def\gnXCIVprim
{  R.Geoghegan, A.Nicas,
\emph{Parameterized Lefschetz-Nielsen fixed
 point theory and Hochshild homology traces},
Amer. J. Math.,
{\bff 116},
(1994)
397--446 }

\def\gnXCV
{  R.Geoghegan, A.Nicas,
\emph{Higher Euler characteristics 1},
L'Enseignement Math\'ematique,
{\bff 41},
(1995)
3--62
}

\def\hadamard
{ J.Hadamard,
\emph{Sur l'it\'eration et les solutions asymptotiques
des \'equations diff\'erentielles},
Bull. Soc. Math. France,
{\bff 29},
(1901)
224 -- 228
}

\def\gokimo
{ H.Goda, T.Kitano, T.Morifuji,
\emph{Reidemeister torsion,
twisted Alexander polynomial
and fibered knots,
}
e-print: math.GT/0311155.
}

\def\gomo
{ H.Goda, T.Morifuji,
\emph{Twisted Alexander polynomial
for $SL(2,\CCC)$-representations
and fibered knots,
} 
C. R. Math. Acad. Sci. Soc. R. Can. 25 (2003) 97-101.
}

\newcommand{\Goda}[1]
{ H.Goda, {#1}}
\newcommand{\Heeg}
{
\emph{Heegaard splitting for sutured
manifolds and Murasufgi sum
},
Osaka J.Math.
{\bff 29},
(1992)
21 -- 40.
}

\newcommand{\handlenumber}
{ 
\emph{On handle number of
Seifert surfaces in $S^3$,
} Osaka J.Math.
{\bff 30},
(1993)
63 -- 80
}

\newcommand{\GodaPaj}[1]
{ H.Goda, A.Pajitnov, {#1}}
\newcommand{\twn}
{
\emph{Twisted Novikov homology 
and circle-valued Morse theory 
for knots and links,} e-print: 
arXiv.math.GT/0312374.
}

\def\hamilton
{ R.S.Hamilton,
\emph{On the inverse
function theorem of Nash and Moser
},
Bull. AMS,
{\bff 7},
(1982)
65 -- 229
}

\def\higman
{ G.Higman,
\emph{Units in group rings},
Proc. London Math. Soc.,
{\bff 46},
(1940)
231 -- 248
}

\newcommand{\huebschmorse}[1]
{W.Huebsch, M.Morse, {#1}}

\newcommand{\bowlbowl}
{
\emph{The bowl theorem and 
a model non-degenerate function,}
Proc. Nat. Acad. Sci U.S.A., 
{\bff 51}, (1964), 49-51.
}

\newcommand{\hulee}
{M.Hutchings, Y-J.Lee
\emph{Circle-valued Morse theory, Reidemeister torsion
and Seiberg-Witten invariants of 3-manifolds},
\quad E-print:
 dg-ga/9612004 3  Dec 1996,
 journal publication:
Topology,
{\bff 38},
(1999),
861 -- 888.
}

\newcommand{\huleee}{ M.Hutchings, Y-J.Lee
\emph{  Circle-valued Morse theory and Reidemeister torsion
},
Geometry and Topology,
{\bff 3},
(1999),
369 -- 396
}

\newcommand{\irwinP}{   M.Irwin,
\emph{ On the Stable Manifold Theorem
},
Bull. London Math. Soc,
{\bff 2},
(1970),
196 -- 198
}

\newcommand{\jiang}{  B.Jiang,
\emph{ Estimation of the
 number of periodic orbits}
Preprint of Universit\"at Heidelberg,
Mathematisches Institut, Heft 65, Mai 1993,
Pac. J. Math. 172, No.1, 151-185 (1996).
}

\newcommand{\kinoshitaterasaka}
{S.Kinoshita, H.Terasaka, 
\emph{ On union of knots},
Osaka math. Journal, 
{\bff 9},
(197), 131 -- 153.}

\newcommand{\kirklivingston}{
P.Kirk, Ch. Livingston,
\emph{
Twisted Alexander invariants,
 Reidemeister torsion, and
 Casson-Gordon invariants}
Topology,
{\bff 38}
1999,
p. 635 -- 661.
}

\newcommand{\kitano}{
T.Kitano,
\emph{
Twisted Alexander polynomial
and  Reidemeister torsion,}
Pacific J.Mat,
{\bff 174,}
1996,
p. 431 -- 442.
}

\newcommand{\kite}{
S.Kinoshita, H.Terasaka,
\emph{
On unions of knots}
Osaka Math. J,
{\bff 9}
1957,
p. 131 -- 153.
}

\newcommand{\klassen}{
E.Klassen,
\emph{
Representations in $SU(2)$
of the fundamental groups
of the Whitehead link
and doubled knots}
Forume Mathematics,
{\bff 5}
1993,
p. 93 -- 109.
}

\newcommand{\kupka}
{   I.Kupka,
\emph{
Contribution \`a la th\'eorie 
des champs g\'en\'eriques},
Contributions to Differential equations,
{\bff 2}
    (1963    ), 457--484,
    {\bff 3}
    (1964    ), 411--420.
 }

\newcommand{\latour}
{F.Latour,
\emph{
Existence de 1-formes ferm\'ees non
singuli\`eres dans une classe de cohomologie 
de de Rham},
Publ. IHES {\bff 80}
(1995),
}

\newcommand{\laudenbach}
{F.Laudenbach,
\emph{On the Thom-Smale complex},~
In {\it An extension of a theorem by Cheeger 
and M\"uller}, ~
Asterisque,{\bff  205} (1992)
p. 219 -- 233.}

\newcommand{\laudsiko}
{   F.Laudenbach, J.-C.Sikorav,
\emph{
Persistance d'intersection avec la section
nulle au cours d'une isotopic hamiltonienne
 dans un fibre cotangent
},
Invent.~Math.~
{\bff 82}
    (1985     ),
pp. 349--357.
 }

\newcommand{\lueckzeta}
{ W.L\"uck,
\emph{The Universal Functorial Lefschetz Invariant }
Fundam. Math. 161, No.1-2, 167-215 (1999).

Preprint: (1998)
}

\newcommand{\liapounov}
{   Liapounov A.M.,
\emph{Probl\`eme g\'en\'eral de 
la stabilit\'e du mouvement},
 Ann. Fac. Sci. Toulouse, 
{\bff 9}, (1907), p. 203-474.
}

\newcommand{\Lin}
{  X.S.Lin,
\emph{
Representations of knot groups
 and twisted Alexander polynomials
},
preprint 1990, 
publication: 
Acta Math. Sin. 
{\bff 17}
    (2001     ),
pp. 361--380.
 }

\newcommand{\McMullen}[1]
{  C.T.McMullen, {#1}}

\def\alexthur
{
\emph{
The Alexander polynomial
of a 3-manifold and the 
Thurston norm on cohomology
},
Ann. Sci. Ec. Norm. Super, 
{\bff 35} No. 2, 153-171 (2002),
http://abel.math.harvard.edu/~ctm/papers.
}

\newcommand{\milncyccov}
{J.Milnor,
\emph{ Infinite cyclic coverings},
In: Conference on the topology of manifolds,
(1968)}

\newcommand{\milnWT}
{J.Milnor,
\emph{ Whitehead Torsion},
Bull. Amer. Math. Soc.
{\bff 72}
(1966),
358 - 426.
}

\newcommand{\peixoto}
{Peixoto,
\emph{ On an approximation theorem of Kupka and Smale},
J. Diff.Eq.
{\bff 3}
(1966),
423 -- 430.
}

\newcommand{\pitcher}
{E.Pitcher,
\emph{ Inequalities of critical point theory},
Bull. Amer. Math. Soc.
{\bff 64}
(1958),
1-30.
}

\newcommand{\noviquasi}
{S.P.Novikov,
\emph{Quasiperiodic Structures in topology},
in the  book:  Topological Methods in Modern Mathematics,
 Publish or Perish,
 1993,
pp. 223--235.  }

\newcommand{\pozniakXC}
{M.Pozniak,
\emph{Triangulation of compact smooth manifolds and Morse theory}
(University of Warwick preprint, 11/1990,
published posthumously as a part of the thesis of M.Pozniak
in Translations of AMS, 2000)}

\newcommand{\pozniakXCI}
{M.Pozniak,
\emph{The Morse complex, Novikov Cohomology 
and Fredholm Theory}
(University of Warwick preprint, 08/1991,
published posthumously as a part
of the thesis of M.Pozniak
in Translations of AMS, 2000)}

\newcommand{\ranicki}[1]
{A.Ranicki, {#1}}

\newcommand{\torsiontwo}
{\emph{The algebraic theory of torsion. I},
Proc. 1983 Rutgers Topology Conference,
Springer Lecture Notes, No. 1126, 199-237}

\newcommand{\torsionone}
{\emph{The algebraic theory of torsion. 
II: products},
K-theory, {\bff 1}, 1987, 115 -- 170}

\def\ranXCV
{ A.Ranicki,
\emph{Finite domination
and  Novikov rings},
Topology,
 {\bff 34} (1995), 619--632.}

\def\ranXCIX
{ A.Ranicki,
\emph{The algebraic construction
of the Novikov complex of a circle-valued
Morse function},

\quad  E-print: math.AT/9903090}

\newcommand{\reidemeister}
{K.Reidemeister,
\emph{Homotopieringe und Linsenr\"aume},
Hamburger Abhandl. {\bff 11} (1938), 102--109.}

\newcommand{\schuetz}[1]
{D.Sch\"utz, {#1}}

\newcommand{\orbs}
{
\emph{Gradient flows of closed 1-forms and
their closed orbits},
e-print:
math.DG/0009055,
~ journal article:
Forum Math. 14(2002) 509--537.
}

\newcommand{\oneparam}
{
\emph{One-parameter fixed point theory 
and gradient flows
of closed 1-forms},
e-print:
math.DG/0104245,
~ journal article:
K-theory, 25(2002), 59-97.
}

\newcommand{\sieben}
{L.Siebenmann,
\emph{A total Whitehead torsion obstruction to 
fibering over the circle},
Comment. Math. Helv. {\bff 45} (1970), 1--48.}

\newcommand{\sikobour}
{J.-Cl; Sikorav,
\emph{
Homologie associ\'ee
a une fonctionnelle
(d'apr\`es Floer)},
Asterisque 201 - 202 - 203, 1991,
S\'eminaire Bourbaki {\bff 733}, 
1990 -- 1991,
November 1991.}

\newcommand{\sikoens}
{ J.-Cl.~Sikorav,
\emph{Un probleme de disjonction par
isotopie symplectique dans un
fibr\'e cotangent},
 Ann.~Scient.~Ecole~Norm.~Sup.,{\bff 19}
 (1986),  543--552.}

\newcommand{\sikoravthese}
{ J.-Cl.~Sikorav,
\emph{
Points fixes de diff\'eomorphismes
symplectiques, intersections de sous-vari\'et\'es
lagrangiennes, et singularit\'es de un-formes ferm\'ees
}
Th\'ese de Doctorat d'Etat Es
Sciences Math\'ematiques,
Universit\'e Paris-Sud, Centre d'Orsay, 1987}

\newcommand{\smale}
{  S.~Smale, \emph{On the structure of manifolds},
 Am.~J.~Math., {\bff 84} (1962)  387--399.}

\newcommand{\smalehcob}
{  S.~Smale, \emph{Generalized Poincar\'e
onjecture in dimensions greater than four},
 Ann.~Math., {\bff 74} (1961)  391--406.}

\newcommand{\smalegraddyn}
{  S.~Smale, \emph{On gradient dynamical systems},
 Am.~J.~Math., {\bff 74} (1961)  199--206.}

\newcommand{\smdyn}
{  S.~Smale,
\emph{Differential dynamical systems},
 Bull. Amer. Math. Soc. {\bff 73} (1967)
747--817.}

\newcommand{\smaletransv}
{  S.~Smale,
\emph{Stable \ma s for differential
equations and diffeomorphisms
},
Ann. Scuola Norm. Superiore Pisa, {\bff 18} (1963)
97--16.}

\newcommand{\smapoi}
{  S.~Smale,
\emph{Generalized Poincare's conjecture in dimensions
greater than four},
 Ann.~Math.,
{\bff 74} (1961)
 391--406.}

\newcommand{\stal}
{J.Stallings ,
\emph{On fibering certain 3-manifolds},
 Proc. 1961 Georgia conference on the Topology of 3-manifolds,
Prentice-Hall, 1962, pp. 95--100.}

\newcommand{\Thang}[1]
{Thang T.Q. Le, {#1}}

\newcommand{\repre}
{
\emph{Varieties of representations and their 
subvarieties of cohomology jumps for knot groups},
(Russian) Mat. Sb. {\bff 184} (1993), p. 57-82.
English translation: Russian Acad. Sci. Sb. Math.
{\bff 78} (1994), p. 187-209.
}

\newcommand{\thom}
{ R.Thom,
\emph{Sur une partition en cellules associ\'ee
\`a une fonction sur une vari\'et\'e},
     Comptes Rendus de l'Acad\'emie de Sciences,
{\bff 228}
(1949),
 973--975.
}

\newcommand{\Thurston}[1]
{W.Thurston, {#1}}

\newcommand{\norm}
{
\emph{A norm for the homology of 3-manifolds},
Mem. Amer. Math. Soc. 
{\bff 59}, (1986), pp. 99 -- 130.
}

\newcommand{\Wada}
{ M.Wada,
\emph{Twisted Alexander polynomial
for finitely presentable groups},
 Topology,
{\bff 33}
(1994),
 241 -- 256.
}

\newcommand{\waldhausen}
{ F.Waldhausen,
 \emph{Algebraic $K$-theory of
 generalized free products, I, II},
Ann. of Math.,
{\bff 108} (1978)
p. 135-204
}

\newcommand{\witten}
{ E.Witten,
 \emph{Supersymmetry and Morse theory},
 Journal of Diff.~Geom.,
{\bff 17} (1985)
 no. 2.
}

\newcommand{\whitehead}
{ J.H.C.Whitehead
 \emph{Simple homotopy types},
 Amer. J. Math.,
{\bff 72} (1952)
 pp. 1- 57
}


\newcommand{\patou}
{ A.V.Pajitnov, \emph{ On the Novikov
complex for rational Morse forms},
\quad preprint:
Institut for Matematik og datalogi, Odense Universitet
Preprints 1991, No 12, Oct. 1991;
~ journal article:
Annales de la Facult\'e de Sciences de
Toulouse {\bf 4}  (1995), 297--338.
The pdf-files available at
http://www.math.sciences.univ-nantes.fr/~pajitnov/
}

\newcommand{\pasur}
{ A.V.Pajitnov,
\emph{
Surgery on the Novikov Complex},
Preprint:
Rapport de Recherche CNRS URA 758,  Nantes,   1993;
the pdf-file available at:
http://www.math.sciences.univ-nantes.fr/~pajitnov/

Journal article:
K-theory {\bff 10} (1996),  323-412.
 }

\newcommand{\pamrl}
{  A.V.Pajitnov,
\emph{Rationality and exponential growth
properties of the boundary operators in the Novikov
Complex},
Mathematical Research Letters,
{\bff 3}
(1996),
  541-548.
 }

\newcommand{\paasym}
{  A.V.Pajitnov,
\emph{   On the asymptotics of
Morse numbers of finite covers of
manifold
},
\quad E-print:
math.DG/9810136, 22 Oct 1998,\quad
journal article:
 Topology,
\textbf{38}, No. 3,  pp. 529 -- 541
(1999).

}

\newcommand{\paadv}
{  A.V.Pajitnov,
\emph{   $C^0$-generic properties of
boundary operators in the Novikov
complex },
\quad E-print:
math.DG/9812157, 29 Dec 1998,
journal article:
Advances in Mathematical Sciences,
 vol. 197, 1999, p.29 -- 117.
}

\newcommand{\pawitt}
{  A.V.Pajitnov,
\emph{Closed orbits of gradient
flows and logarithms of non-abelian Witt vectors},
\quad E-print:
 math.DG/9908010, 2 Aug. 1999
journal article:
  K-theory, Vol. 21 No. 4, 2000.
}

\newcommand{\pajandran}
{A.V.Pajitnov, A.Ranicki,
\emph{The Whitehead group of the Novikov ring},
\quad E-print:
 math.AT/0012031, 5 dec 2000,
 journal article:
  K-theory, Vol. 21 No. 4, 2000.
}

\newcommand{\pagrad}
{A.V.Pajitnov,
\emph{$C^0$-topology in Morse theory},
\quad E-print:
 math.DG/0303195, 5 dec 2000.
}

\newcommand{\godapaj}
{H.Goda, A.V.Pajitnov,
\emph{Novikov homology and fibering obstructions
for knots and links},
in preparation
}

\newcommand{\paodense}
{ A.V.Pazhitnov,
\emph{
On the Novikov
complex for rational Morse forms
}, Preprint:
~
Institut for Matematik og datalogi
 Odense Universitet,
 Preprints 1991, No. 12
 Odense, October 1991.
 }


\newcommand{\postnTM}
{   M.M.Postnikov,
\emph{Introduction to Morse theory
},
 Moscow, Nauka,  1971, in Russian
}


\newcommand{\arnoldEquadiff}
{   V.I.Arnold,
\emph{Ordinary Differential Equations
},
 Moscow, Nauka,  1975.
}

\newcommand{\abrob}
{   R. Abraham, J.Robbin,
\emph{Transversal mappings and flows
},
 Benjamin, New York,  1967.
}

\newcommand{\atiyahandmacdo}
{   M.F.Atiyah, I.G.Macdonald
\emph{Introduction to commutative algebra
},
Addison-Wesley,   1969.
}

\newcommand{\bass}
{ H.Bass,
\emph{Algebraic K-theory},
Benjamin, 1968.
}

\newcommand{\birkrota}
{ G.Birkhoff, G-C. Rota,
\emph{Ordinary differential equations},
Blaisdell Publishing Company, 1962.
}

\newcommand{\bz}
{ G.Burde, H.Zieschang,
\emph{Knots},
de Gruyter, 1985.
}

\newcommand{\bour}
{N.Bourbaki,
\emph{Groupes de Lie, Alg\`ebres de Lie
}
}

\newcommand{\cartanCD}
{   H.Cartan
\emph{Cours de Calcul Diff\'erentiel
},
Hermann,  1977.
}

\newcommand{\cohen}{
M.M.Cohen,
\emph{A course in Simple-Homotopy theory},
Springer, 1972.}

\newcommand{\crfox}{
M.H.Crowell, R.H.Fox,
\emph{Introduction to knot theory},
Ginn and company, 1963.}

\newcommand{\dieudonne}
{   J.Dieudonn\'e
\emph{Foundations of modern analysis
},
Academic press, 1960.
}

\newcommand{\dieud}
{   J.Dieudonn\'e
\emph{El\'ements d'analyse, Vol. III
},
Gauthier-Villars, 1970
}

\newcommand{\dold}
{   A.Dold
\emph{Lectures on Algebraic Topology
},
Springer,  1972.
}

\newcommand{\eisenbud}
{   D.Eisenbud
\emph{Commutative Algebrs with a View Toward
 Algebraic Geometry
},
Springer,  1994.
}

\newcommand{\irwinB}
{ M.C.Irwin
\emph{Smooth dynamical systems
},
Academic Press,  1980.
}

\newcommand{\hartman}
{ P.Hartman
\emph{Ordinary differential equations},
Wiley; New York,  1964.
}

\newcommand{\hirsch}
{ M.Hirsch
\emph{Differential Topology
},
Springer,  1976.
}

\newcommand{\huse}
{ D.Husemoller
\emph{Fibre bundles
},
McGraw-Hill, 1966.
}

\newcommand{\kelley}
{ J.L.Kelley
\emph{General Topology
},
van NOstrand, 1957.
}

\newcommand{\kiang}
{   Kiang Tsai-han
\emph{The theory of Fixed point classes},
 Springer, 1989.
}

\newcommand{\kling}
{   W.~Klingenberg
\emph{Lectures on closed geodesics
},
 Springer, 1978.
}

\newcommand{\lafont}
{  S.Gallot, D.Hulin, J.Lafontaine,
\emph{Riemannian Geometry },
 Springer.
}

\newcommand{\McDonald}[1]
{B.R.McDonald, {#1}}

\newcommand{\linalg}
{
\emph{Linear Algebra Over Commutative Rings,}
Marcel Dekker, 1984, 544 pages.}

\newcommand{\milnMT}
{   J.~Milnor,
\emph{ Morse theory
},
 Princeton University Press, 1963.
}

\newcommand{\milnKT}
{   J.~Milnor,
\emph{ Introduction to algebraic K-theory
},
 Princeton University Press, 1971.
}

\newcommand{\milnhcob}
{   J.~Milnor,
\emph{Lectures on the
h-cobordism theorem},
 Princeton University
Press,
 1965.
}

\newcommand{\milnstash}
{   J.Milnor and J.Stasheff,
\emph{ Characteristic Classes},
 Princeton University Press,
 1974.}

\newcommand{\morse}[1]
{M.Morse, {#1}}

\newcommand{\calcvar}
{
\emph{Calculus of Variations in the Large},
  American Mathematical Society Colloquium 
Publications,
Vol.18,
 1934.
}

\newcommand{\munkres}
{  J.R.Munkres,
\emph{Elementary differential toplogy},
Annals of Math. Studies,
Vol.54, Pinceton
 1963.}

\newcommand{\cohn}
{ P.M.Cohn,
\emph{Free rings and their relations},
  Academic press
( 1971)}

\newcommand{\lam}
{    T.Y.Lam,
\emph{Serre's Conjecture,         }
Lecture Notes in Mathematics {\bff 635}, (1978) 227 p.
}

\newcommand{\lang}
{    S.Lang,
\emph{Algebra ,         }
Addison-Wesley (1965)
}

\newcommand{\massey}
{   W.Massey,
\emph{ Homology and cohomology theory
},
 Marcel Dekker, 1978.
}

\newcommand{\palisdemelo}
{   J.Palis, Jr.,
W.de Melo,
\emph{ Geometric  theory of dynamical systems},
 Springer, 1982.
}

\newcommand{\ranKL}
{   A.A.Ranicki,
\emph{Lower $K$- and $L$-theory,         }
LMS Lecture Notes 178, Cambridge, 1992
}

\newcommand{\ranKNO}
{   A.A.Ranicki,
\emph{High-dimensional knot theory,         }
Springer, 1998
}

\newcommand{\ranitor}
{   A.A.Ranicki,
\emph{The algebraic theory of torsion I.,}
Lecture Notes in Mathematics {\bff 1126} (1985), 199--237.
}

\newcommand{\rock}
{   R.T.Rockafellar,
\emph{Convex Analysis},
Princeton University Press (1970).
}

\newcommand{\rolfsen}
{   D.Rolfsen,
\emph{Knots and Links},
Publish or Perish  (1976, 1990).
}

\newcommand{\rosenberg}
{   J.Rosenberg,
\emph{Algebraic $K$-theory
 and its applications},
Springer, (1994).
}

\def\sharko{ V.V.Sharko,
\emph{ Funktsii na mnogoobraziyah
(algebraicheskie i topologicheskie aspekty)}, Kiev,
Naukova dumka
(1990),
 31-35.
   \quad English translation:
 }

\newcommand{\spanier}
{ E.H.Spanier,
\emph{Algebraic topology},
  McGraw-Hill,
(1966)}

\newcommand{\stong}
{R.Stong,
\emph{Notes on Cobordism theory},
Princeton, New Jersey, 1968}

\newcommand{\switzer}
{ R.M.Switzer,
\emph{Algebraic topology -- homotopy and homology},
  Springer,
( 1975)}

\newcommand{\wander}
{ B.L.van der Waerden,
\emph{Algebra 1},
Springer, 
( 1971)}

\newcommand{\wells}
{ R.O.Wells,
\emph{Differential Analysis on complex manifolds},
Prentice-Hall,
( 1973)}

\newcommand{\weibel}
{ C.A.Weibel,
\emph{An introduction to homological algebra },
Cambridge University press,
( 1997)}



\newcommand{\burghelea}{ D.Burghelea, S.Haller
\emph{ On the topology and analysis of a closed 1-form
(Novikov's theory revisited)
}, \quad  E-print dg-ga/0101043  5 jan. 2001}

\newcommand{\hulihuli}{ M.Hutchings, Y.J.Lee,
\emph{ Circle-valued Morse theory and Reidemeister Torsion
}, \quad  E-print dg-ga/9706012 23 June 1997}

\newcommand{\mengt}
{ G.Meng, C.H.Taubes
\emph{SW=Milnor Torsion }
\quad   preprint
(1996)
}

\newcommand{\ranprepr}
{  A.A.Ranicki,
      \emph{
Finite domination and Novikov rings},
preprint,
 1993 }

\newcommand{\sheiham}
{  D.Sheiham,
\emph{Noncommutative characteristic
polynomials and Cohn localization},
\quad   preprint, 2000}

\newcommand{\minervThesis}
{ G.Minervini
\emph{ A current approach to Morse and Novikov theories}
Tesi di Dottorato, Universit\`a "La Sapienza", Roma}




\def\farber{M.Farber
 \emph{
 Exactness of Novikov inequalities },
Functional Analysis and Applications  {\bff
 19},
 1985.
 }

\def\kozlovskii{O.S.Kozlovskii
\emph{The dynamics
of intersections of
analytic manifolds
}, Doklady ANSSSR,
{\bf 323}
(1992).
\quad English translation:

 Sov.Math.Dokl.
{\bff 45}
(1992), 425--427.}

\def\novidok{S.P.Novikov,
\emph{ Many-valued functions
 and functionals. An analogue of Morse theory  },
Doklady AN SSSR,
{\bf 260}
(1981),  31-35 (in Russian),
\quad English translation:
 Sov.Math.Dokl.
{\bff 24}
(1981), 222-226. }

\newcommand{\noviuspe}
{ S.P.Novikov,
\emph{The hamiltonian formalism and a
multivalued analogue of
Morse theory,
}
Uspekhi Mat. Nauk,
{\bff 37}    (1982),  3-49(in Russian),
\quad English translation:
 Russ. Math. Surveys,
{\bff 37} (1982),
 1--56.}

\newcommand{\pasbor}
{  A.Pajitnov, 
\emph{
O tochnosti neravenstv tipa Novikov dlya mnogoobrazii
so svobodnnoi abelevoi fundamental'noi gruppoi},
Mat. Sbornik, 
(1989)
no. 11.

\quad English translation:

 A.V.Pazhitnov,\emph{
 On the sharpness of
 Novikov-type inequalities for
manifolds with free abelian fundamental group.},
 Math. USSR Sbornik,
{\bff 68}
 (1991),
 351 - 389.
}

\def\pazamet{
A.V.Pazhitnov,
\emph{Modules over some localizations
of the ring of Laurent polynomials},
Mathematical Notes, {\bff 46}, 1989, no. 5.}

\newcommand{\pafest}
{  A.V.Pajitnov,
\emph{ Simple homotopy type of Novikov complex
and $\zeta$-function of the gradient flow, }
\quad E-print:
dg-ga/970614 9 July 1997;
journal article:
Russian Mathematical Surveys,
\textbf{54}
(1999), 117 -- 170.}

\newcommand{
\pastpet}
{  A.V.Pajitnov,
{\it   The incidence coefficients in the Novikov
complex are
generically rational functions,}
\quad E-print:  dg-ga/9603006 14 March  96,
journal article: Algebra i Analiz,
{\bff 9}, no.5 (1997),  92--139 (in Russian),
\quad English translation:
Sankt-Petersbourg Mathematical Journal
\textbf{9}
(1998),
no. 5, p. 969 -- 1006.
}

\newcommand{\prw}
{  A.V.Pajitnov, C.Weber, L.Rudolph,
\emph{ Morse-Novikov number for knots and links},~
Algebra i Analiz,
{\bff 13}, no.3 (2001),
(in Russian),
English translation:
Sankt-Petersbourg Mathematical Journal.
{\bff 13}, no.3 (2002), p. 417 -- 426.
}

\newcommand{\paclo}
{  A.V.Pajitnov,
\emph{ Counting closed orbits
of gradients of circle-valued maps},~
E-print:  math.DG/0104273 28 Apr. 2001,
journal article:
Algebra i Analiz,
{\bff 14}, no.3 (2002),  92--139
(in Russian),
English translation:
Sankt-Petersbourg Mathematical Journal.
{\bff 14}, no.3 (2003).
}

\newcommand{\turaev}[1]
{V.Turaev, {#1}}

\def\alex
{
\emph{ The Alexander polynomial 
of a three-dimensional manifold},
Math. USSR Sbornik, v. 26 (1975), No. 3, p. 313-329.
}

\def\reidknot
{
\emph{
Reidemeister Torsion in knot theory,}
Russian Math. Surveys, 41:1 (1986), 119 - 182.
}

\newcommand{\euler}
{
\emph{ Euler structures, nonsingular vector fields
and torsions of Reidemeister type},
Math. USSR Izvestia 34:3 (1990), 627 - 662.
}

\def\twocomple
{
\emph{ A norm for the cohomology of 2-complexes},
Algebraic and Geometric topology, 
Vol. 2(2002), 137-155.
}

\title[Novikov homology and  Thurston cones]
{Novikov homology, 
twisted Alexander polynomials
and Thurston cones\\
}
\author{A.V.Pajitnov}
\address{Laboratoire Math\'ematiques Jean Leray 
UMR 6629,
Universit\'e de Nantes,
Facult\'e des Sciences,
2, rue de la Houssini\`ere,
44072, Nantes, Cedex}
\email{ pajitnov@math.univ-nantes.fr}
\begin{abstract}
Let $M$ be a connected CW complex,
let $G$ denote the fundamental group
of $M$. Let $\pi$
be an epimorphism of $G$ onto a free
finitely generated abelian group $H$, let
  $\xi:H\to \mathbf R$ be 
a homomorphism and  $\r$ be an anti-homomorphism
of $G$ to the group $GL(V)$ of automorphisms of a free
finitely generated $R$-module $V$
(where $R$ is a commutative
factorial ring). 

We associate to these data the 
{\it twisted Novikov homology}
of $M$ which is a module over the Novikov  
completion of the ring $\L=R[H]$. 
The twisted Novikov homology 
provides the lower bounds for the number of zeros of 
any Morse form which cohomology class equals 
$\xi\circ\pi$.
This construction generalizes  the
work by H.Goda and the author (arXiv:math.DG/0312374).

In the case when $M$ 
is  a compact connected 
3-manifold $M$ with zero Euler characteristic
we obtain a criterion of vanishing of the 
twisted Novikov homology of $M$
in terms of the corresponding 
twisted Alexander polynomial
of the group $G$.

We discuss the relations of  the 
twisted Novikov homology with the Thurston norm on 
the 1-cohomology of $M$.
\end{abstract}

\maketitle


\section{Introduction}
\mlb{s:intro}

Let $M$ be a closed manifold,
 $f:M\to S^1$
be a circle-valued  Morse function
on $M$.
Let $m_k(f)$ denote the number 
of critical points of $f$
of index $k$.
The Morse-Novikov theory 
provides lower bounds for the 
numbers $m_k(f)$
which are computable in terms of the 
homotopy type 
of $M$ and the homotopy class of $f$.
The general schema of obtaining such 
bounds is as follows
(see \cite{novidok}, \cite{farber}, \cite{pasbor}, 
\cite{patou}).
Consider a regular covering $\PP:\bar M\to M$
with structure group $G$, such that 
the function 
$f\circ\PP:\bar M\to S^1$
is homotopic to zero (or equivalently, $f$
lifts to a Morse function 
$\bar M\to \RRR$).
The induced \ho~ $f_*:\pi_1(M)\to\ZZZ$
can be  factored through a \ho~
$\xi=\xi(f):G\to\ZZZ$.

Applying the standard method of counting 
gradient
flow lines (see \cite{witten})
one  obtains a chain complex $\NN_*$
({\it the Novikov complex})
over a certain completion $\LLLxi$ 
of the group ring $\LL=\ZZZ[G]$.
\footnote{See 
the definition of  $\LLLxi$ in 
Subsection \mrf{su:nov_hom}.})
The Novikov  complex is freely 
generated over $\LLLxi$
by the critical points
of $\o$, and its homology  is 
isomorphic  to  the homology of  
the tensor product 
\bq \lb{f:sing_comple}
\wh \CC_*(\bar M, \xi)
=
\LLLxi\tens{\LL} \CC_*(\bar M),
\end{equation}
where $\CC_*(\bar M)$ is
the cellular  chain complex 
of the  covering $\bar M$.
In particular, if 
the chain complex \rrf{f:sing_comple}
is not contractible, 
the function $f$
must have at least one critical point.
Developing further this observation,
one can obtain 
lower bounds for the numbers $m_k(f)$
in terms of the 
numerical invariants of the 
homology of the chain complex 
\rrf{f:sing_comple}.

Let us discuss different possible choices 
of the covering $\PP$.
The universal covering $\wi M\to M$ 
contains certainly the 
maximum of information.
The disadvantage is that the corresponding
ring $\LLLxi$, being a completion of the group ring 
of the fundamental group 
may be  very complicated, and it can be  
difficult to extract 
the necessary numerical data.

Another obvious possibility is
the infinite cyclic covering 
$$\PP_f:\bar M_f\to M$$
induced by $f$ from the universal 
covering $\RRR\to S^1$.
Here the ring $\LLLxi$
is a principal ideal domain, 
and the explicit lower bounds 
for the numbers of the critical points 
are  easy to deduce
(see Subsection \mrf{su:nov_hom}
for more details).

An intermediate choice 
is the maximal free abelian covering
$$\PP_{ab}:\ove M_{ab}\to M.$$
The structure group of this covering 
is equal to $H_1(M)/Tors$.
The Novikov ring in this case is a 
completion of the Laurent polynomial
ring in several variables, and its 
homology properties are in general
rather complicated. 
But this choice has an advantage
that for any function $f:M\to S^1$ 
the function $f\circ \PP_{ab}$ 
is homotopic to zero.
This allows to study the dependance of 
the Novikov homology
on the class 
$$\xi=\xi(f)
\in \Hom(H_1(M)/Tors, \ZZZ)\approx H^1(M,\ZZZ).$$
In particular one can 
get some information about the set of all
$\xi$
such that the Novikov homology
$$
\wh H^{ab}_*(M,\xi)
=
H_*\Big(\wh\CC_*(\ove M_{ab}, \xi)\Big)
$$ 
vanishes
and, therefore, obtain an information 
about the set of
classes in $H^1(M,\ZZZ)=[M,S^1]$
representable by fibrations.

\bede\mlb{d:integral_c}
Let $H$ be a finitely 
generated free abelian group.
Put
$$
H_\RRR=H\otimes \RRR,
\quad
H'_\RRR=\Hom(H,\RRR)
=
\Hom(H,\ZZZ)\otimes \RRR.
$$
A {\it closed  cone} in $H'_\RRR$
is a closed subset $C$ such that
$v\in C\Rightarrow \l \cdot v\in C$
for every $\l\geq 0$.

An {\it open cone} in $H'_\RRR$
is an open subset $D$ such that
$v\in C\Rightarrow \l\cdot v\in C$
for every $\l > 0$.

An {\it integral hyperplane } of
the vector space 
$H'_\RRR=\Hom(H,\RRR)$
is a vector subspace of codimension 1 
in $H'_\RRR$ having a basis formed by
elements of $\Hom(H, \ZZZ)$.

A connected component of 
the complement to a given hyperplane $\G$
will be called {\it open half-space}
corresponding to  $\G$. The closure 
of an open 
half-space will be called
{\it closed half-space} corresponding 
to  $\G$.

A closed cone which 
is the intersection of a finite family of closed 
half-spaces
corresponding to integral hyperplanes
is called  {\it  closed polyhedral cone}.

An open  cone which 
is the intersection of a finite family 
of open half-spaces
corresponding to integral hyperplanes
is called  {\it  open  polyhedral cone}.

A subset $C\sbs H'$
is called 
{\it open  polyhedral conical subset, }
if it is empty or equals 
$H'\sm \{0\}$
or is a finite disjoint union 
of open polyhedral cones.

A subset $A\sbs H'$
is called  {\it  
quasi-polyhedral conical subset }
if there is an 
open  polyhedral conical subset $C$
and a finite union $D$
of integral hyperplanes,
such that $D\cup C=D\cup A$.

\end{defi}

The next theorem follows from the main theorem of 
my paper  \cite{pazamet},
it is based on  earlier 
results of J.-Cl. Sikorav
(see \cite{pasbor}).
\beth\mlb{t:cones_nontwi}
The set of all 
classes $\xi\in H^1(M,\ZZZ)$
such that the Novikov homology
$\wh H^{ab}_*(M,\xi)$
vanishes
is the intersection with $H^1(M,\ZZZ)$
of a quasi-polyhedral conical subset of
$H^1(M,\RRR)$.
\footnote{ The results  of \cite{pasbor} and
\cite{pazamet} pertain actually to a 
more general case of 
arbitrary homomorphisms $\pi_1(M)\to\RRR$,
and not only \ho s $\pi_1(M)\to\ZZZ$, 
see the discussion below. })
\enth

Now let us proceed to 
non-abelian  coverings.
In the joint work with H.Goda \cite{GodaPaj}
we introduced a new version of the Novikov 
homology. We call it {\it twisted Novikov homology}.
The input data for the construction
is: a connected CW 
complex
 $M$, a \ho~  $\xi:\pi_1(M)\to\ZZZ$,
and an anti-homomorphism
$\rho:\pi_1(M)\to GL(n,\ZZZ)$. The resultant
twisted Novikov homology groups are  
modules  over the 
principal ideal domain   $\ZZZ((t))$,
so the numerical invariants 
are easily extracted from
the homological data. 
On the other hand 
the non-abelian homological algebra
of the universal covering of $M$
is encoded in this homology via 
the representation $\rho$.
The construction of the twisted Novikov homology
is motivated by the notion of 
{\it twisted Alexander polynomial}
for knots and links.  (See the papers 
\cite{Lin} of X.S.Lin and  \cite{Wada} 
of M.Wada  for the definition and properties of
the twisted Alexander polynomials, 
and the paper \cite{gokimo} of H.Goda,
T.Kitano, T.Morifuji for applications of 
\talp s to fibering 
obstructions for knots and links).

The definition of the twisted Novikov 
homology generalizes immediately 
to the case of arbitrary cohomology classes 
$\xi\in H^1(M,\RRR)$.
The input data for this construction 
is as follows. Let $R$ be a commutative ring.
Let $V$ be a   \fg~free left $R$-module.
Denote by $GL_R(V)$ the group of all
automorphisms of $V$ over $R$.
Let $\r:G\to GL_R(V)$
be an anti-homomorphism
(that is, $\r(ab)=\r(b)\r(a)$
for all $a,b\in G$; $\r$ will also be called
{\it right representation}). 
Let $\pi:G\to H$
be an epimorphism of $G$ onto
a free finitely generated abelian group $H$.
Let $\xi:H\to\RRR$
be a group \ho.

To this data we 
associate the twisted 
Novikov homology as follows.
Let $\L=R[H]$; put $V^H=\L\tens{R}V$;
define a right representation
$\r_\pi:G\to GL_\L(V^H)$ as follows:
\bq\lb{f:r_pi}
\r_\pi(g) \Big(\sum_i\l_i\otimes v_i\Big)
=
\sum_i \big(\pi(g)\l_i\big)\otimes \r(g)v_i.
\end{equation}
Form the tensor product of the cellular chain complex
$\CC_*(\wi M)$ of the universal covering
with the right $\ZZZ G$-module $V^H$:
$$
\wi \CC_*(M,\r_\pi)
=
V^H\tens{\ZZZ G} \CC_*(\wi M).
$$
This is a chain complex of left free 
$\L$-modules (observe that $\rk_\L(\wi \CC_k(M,\r_\pi))=
n\cdot \rk \CC_k(\wi M)$,
where $n$ is the rank of $V$ over $R$).
Apply the tensor product with $\Lxi$
to obtain the chain complex
\bq\lb{f:nov_twi_com}
\wh \CC_*(M,\r_\pi,\xi)
= \Lxi\tens{\L}
\wi \CC_*(M,\r_\pi).
\end{equation}

Its homology 
\bq\lb{f:h_nov_twi_com}
\wh H_*(M,\r_\pi,\xi)=
H_*(\wh \CC_*(M,\r_\pi,\xi))
\end{equation}
is called {\it twisted 
Novikov homology}.
The twisted Novikov homology of
our paper \cite{GodaPaj}
corresponds to the particular case 
when  $H=\ZZZ$, and the homomorphism 
$\xi:H\to \RRR$ above is the inclusion 
$\ZZZ\rInto \RRR$.

The present  generalization makes simpler the 
statements of several  theorems below.
It  has also a geometrical background; 
it corresponds to  Morse forms
while the framework of our paper \cite{GodaPaj}
was related to circle-valued Morse functions.
(Recall that a closed 1-form on a manifold
$M$ is called {\it Morse form},
if locally it is a differential 
of a Morse function.)
Namely we have the following theorem:
\beth\mlb{t:inequa}
Let $\o$ be a Morse form on a closed 
connected manifold $M$.
 Assume that 
the cohomology class
$[\o]\in H^1(M,\RRR)=Hom(\pi_1(M), \RRR)$
can be factored as
$[\o]=\xi\circ\pi$,
where $\pi:\pi_1(M)\to H$
is an epimorphism onto 
a free abelian group,
and $\xi:H\to\RRR$ is a homomorphism.
Then there is a chain complex 
$\NN^\r_*$ such that
\belis\item
$\NN^\r_k$
is a free $\Lxi$-module 
with $n\cdot m_k(\o)$ 
free generators in degree $k$
(where $m_k(\o)$ stands for the number
of zeros of the form $\o$ of index $k$).
\item
$\NN^\r_*$
is chain homotopy equivalent to 
$\wh \CC_*(\wi M,\r_\pi,\xi)$.
\enlis
\enth
In particular, if the cohomology 
class $[\o]$
contains a nowhere vanishing 1-form, then 
the twisted Novikov homology 
$\wh H_*(M,\r_\pi, \xi )$
equals to zero. 

The most natural choice of
the epimorphism $\pi:\pi_1(M)\to H$
is the projection $\pi_1(M)\to H_1(M)/tors$
onto the integral homology group of $M$ 
modulo its torsion subgroup.
The corresponding twisted Novikov 
 homology will be denoted by 
$\wh H_*(M, \r,\xi)$.

\bede\mlb{d:acyc}
A non-zero homomorphism 
$$\xi\in Hom(H_1(M),\RRR)=
\Hom(H_1(M)/Tors, \RRR) =
H^1(M,\RRR)$$
is called {\it $\r$-acyclic}, if the  
$\r$-twisted Novikov homology 
$\wh H_*(M,\r, \xi)$
vanishes.
The set of all $\r$-acyclic 
classes will be denoted
$\VV_{alg}(M,\r)$.
\end{defi}
The reason for 
studying the $\r$-acyclic classes
is that any class $\xi$ containing 
a nowhere vanishing 
closed 1-form is $\r$-acyclic 
for any representation $\r$.
An immediate generalization of Theorem 
\mrf{t:cones_nontwi}
leads to the following result about the 
algebraic structure
of the set of $\r$-acyclic classes.
\beth\mlb{t:cones_twi}
For a given 
right representation $\r$ of $\pi_1(M)$
the set $\VV_{alg}(M,\r)$ 
of all $\r$-acyclic classes $\xi$
is a  quasi-polyhedral conical subset.
\enth

In general we do not know whether 
the set of all $\r$-acyclic classes is an
open polyhedral conical subset.
In other words, we can describe the structure
of the set of all $\r$-acyclic classes only up
to some finite union of 
hyperplanes in $H^1(M,\RRR)$.
However  in the case when $M$ is a 
3-manifold, 
we have 
a much stronger assertion, 
which is the main result of 
the present paper (see Section \mrf{s:three}):
\beth\mlb{t:cones_twi_3}
Let $M$ be a  connected compact 
three-dimensional
manifold (maybe with a non-empty boundary),
such that $\chi(M)=0$.
Let $\r$ be a right representation of $\pi_1(M)$.
Then:
\beli\item 
The set $\VV_{alg}(M,\r)$ 
of all $\r$-acyclic classes $\xi$
is an open  polyhedral conical subset.
\item
This subset is entirely determined by the \talp~
associated to the group $\pi_1(M)$
and the representation $\r$.
\enli
\enth

The open polyhedral cones
forming the subset $\VV_{alg}(M,\r)$ 
 will be called
{\it $\r$-acyclicity cones}.
The theorem above implies that 
the set of all $\r$-acyclic classes 
depends only on the group $\pi_1(M)$
and the representation $\r$.

Along with the \talp s
the proof of the above theorem uses 
another polynomial invariants 
of the chain complexes,
which we introduce in Section
\mrf{s:fitt_invar}
and call {\it the Fitting  invariants}.
These invariants are defined as the GCDs of 
the minors 
of the second boundary operator of the chain 
complex 
\rrf{f:nov_twi_com}.
We show that they 
are directly related 
to the Novikov
homology, and in many cases the Novikov 
homology in degree one can be computed 
from the sequence of the  Fitting invariants.
The \talp~is a much more sophisticated invariant,
but it turns out that the 
image of the \talp~in the Novikov completion $\Lxi$
is essentially the same as the image of the 
corresponding Fitting invariant.

Theorem \mrf{t:cones_twi_3}
is related to the famous Thurston's
theorem \cite{th}, which implies  
that the set $\VV(M)$ of all 
classes $\xi\in H^1(M,\RRR)$
representable by closed nowhere  vanishing 
1-forms is a finite union of open 
polyhedral cones, namely the cones
on certain faces of the unit ball 
of the Thurston norm on $H^1(M,\RRR)$.
We shall call these cones {\it Thurston cones}.
For every right 
representation $\rho$
of $\pi_1(M)$ in $GL(\ZZZ^n)$
the set $\VV(M)$ is contained in 
the set
$\VV_{alg}(M,\r)$ of all $\r$-acyclic classes:  
$$
\VV(M)\sbs \VV_{alg}(M,\rho)
$$
so each of the  Thurston  cones 
is  contained in one of the 
$\r$-acyclicity cones.
The $\r$-acyclicity cones 
 are    computable
from the twisted Alexander polynomial,
which is in its turn
computable from the Alexander matrix 
associated to any finite presentation 
of the group $\pi_1(M)$.
Thus the set $\VV_{alg}(M,\rho)$
is in a sense a computable upper bound for
the set $\VV(M)$. Let 
$$\VV_{alg}(M)
=
\bigcap_{\r}\VV_{alg}(M,\r)
$$
where $\r$ ranges over the set 
of all right representations 
$\r$
of $\pi_1(M)$ in $GL(\ZZZ^n).$
We have then:
\bq\lb{f:geo_alg}
\VV(M)
\sbs
\VV_{alg}(M).
\end{equation}
It is an interesting problem to 
investigate the relationship between the two sets,
and in particular answer the next question:
\pa
{\it 
For which manifolds 
the equality 
$\VV(M) = \VV_{alg}(M)$
holds? }
\pa
We do not know examples of manifolds
for which $\VV(M)\not= \VV_{alg}(M)$.
On the other hand, it is easy to construct 
manifolds for which 
$$
\VV(M)
\not=
\bigcap_{\r\in \RR_f}\VV_{alg}(M,\r),
$$
where $\RR_f$
is the set of all right representations
of $G$ over finite fields
(see Section \mrf{su:nov_thur}).
For such  manifolds the
right representations over finite fields are
not sufficient to
detect all the cohomology classes 
representable by non-singular 1-forms

A certain amount of computations
will be necessary to clarify 
the relation between $\VV(M)$ and $\VV_{alg}(M)$
and  answer  the question above.
The recent progress in the software related
to the computations of invariants of  
knots and links, especially 
the Kodama's KNOT program allows 
us to hope that such computations 
can be carried out.

\centerline{\bf Acknowledgements.}
A part of the
work presented in this paper was done  during 
the author's stay
at Osaka City University in December 2003.
I was motivated  by many discussions
with Professors Kunio Murasugi, 
Akio Kawauchi, Makoto Sakuma, Masaaki Wada,
Hiroshi Goda, Teruaki Kitano, 
Kouji Kodama, Takayuki Morifuji,
Mikami Hirasawa. 
In particular, the question 
about the  relation between 
the Novikov homology and the twisted Alexander
polynomial was raised during these discussions. 
The present paper contains a partial answer 
to this and other questions of algebraic
nature related to the twisted Novikov homology.
I am grateful to Osaka City University 
and to Professor Akio Kawauchi 
for the warm hospitality during my stay.

The present paper is a second part 
of  our joint work 
\cite{GodaPaj} with H.Goda. 
My special thanks  to Hiroshi Goda
for sharing his knowledge and insight 
in knot theory.

Thanks to A.Ranicki and Liam O'Carroll
for providing useful references on the Fitting
invariants of modules.
\pa

\centerline{\bf Notes on the terminology.}

{\it Ring} means an associative ring with a unit.
For a ring $R$ we denote by $R^\bu$ 
the multiplicative group
of all invertible elements of $R$.
{\it Module} means always a left module,
if the contrary is not
stated explicitly.
The homology of a space $X$ with 
integral coefficients
will be denoted
by $H_*(X)$.
For a left module $V$ over a commutative ring $R$
we denote by $GL_R(V)$
the group of all $R$-automorphisms of $V$.
When the ring $R$ is clear from the context,
we use the abbreviated notation $GL(V)$.
For a ring $R$ the symbol $GL(n,R)$
denotes the multiplicative group of
all invertible $n\times n$-matrices 
with coefficients in $R$.


\section{Twisted Novikov homology}
\mlb{s:nov_twi}
${}$

\subsection{ Novikov homology}
\mlb{su:nov_hom}

\bede\mlb{d:nov_ring}
Let  $G$ be a  group, and $R$ be a 
commutative ring with a unit.
Put $\LL=R[G]$.
Let $\xi:G\to\RRR$ be a group homomorphism.
Let $\wh{\wh{ \LL}} $
be the set of all formal linear 
combinations (infinite 
in general) of the form
$\l=\sum_{g\in G} n_g g$
where $n_g\in R$.
For $\l\in \wh{\wh\LL}$ and $C\in\RRR$
put
$$\supp(\l, C)=\{g\in G~|
~n_g\not=0 \quad \& \quad \xi(g) > C\}.$$
The Novikov ring is defined as follows:
$$
\LLLxi=
\{\l\in  \wh{\wh\LL} ~|~ 
\supp(\l, C) 
\mxx{ is finite
for every }  C\in\RRR\}.
$$
(It is easy to show that the subset
$\LLLxi\sbs \wh{\wh\LL}$
has indeed a natural structure of a ring, 
containing $\LL$ as a subring.)
\end{defi}

\begin{exam}\mlb{rat}
Let $G$ be an abelian finitely generated group.
A homomorphism $\xi:G\to\RRR$
is called {\it rational}
if $\xi=\l\cdot \xi_0$, where 
$\xi_0:G\to\ZZZ$ is a \ho, and $\l\in\RRR$.
Equivalently, $\xi$ is rational, if $\xi=0$
or $\xi(G)\approx \ZZZ$.
When $\xi$ is non-zero and 
rational, we have an isomorphism
$$
\LLLxi\approx K[[t]][t^{-1}]
\mxx{ where } K=R[\Ker \xi]
$$
so that the ring $\LLLxi$
is a localization of the ring of the power series
$K[[t]]$.
In particular $\LLLxi$ is Noetherian.
\end{exam}

\begin{exam}
\mlb{e:tot_irrat}
Let $G$ be a free abelian finitely generated group.
In this case 
$\LL=R[G]$ is isomorphic to the ring of Laurent 
polynomials in $k$ variables 
(where $k=\rk G$)
with  coefficients in $R$.

When ~$\xi:G\to \RRR$
is  
 injective (such \ho s are also called 
{\it totally irrational}),
the algebraic 
properties of the ring $\LLLxi$
are suprisingly simple: 
\beth\mlb{t:siko1}
If $\xi$ is totally irrational, 
and $R$ is a principal ideal domain, 
then the ring $\LLLxi$ 
is also a principal
ideal domain.
\enth

For the case $R=\ZZZ$ 
this  theorem is  due to J.-Cl. Sikorav
(his proof was published in \cite{pasbor}).
The proof in the general case is similar.
Observe that if $R$ is a field, then 
$\LLLxi$ is also a field.

\end{exam}

The main topological applications of these 
constructions are in
the theory of Morse forms.
Let $M$ be a closed connected $\smo$ manifold,
let $\o$ be a Morse form on $M$.
Let $[\o]\in H^1(M,\RRR)$
denote the de Rham cohomology class 
of $\o$;
then $[\o]$ can also be identified with a 
homomorphism
$G=\pi_1(M)\to\RRR$.
Let $\wh\gL_{[\o]}$
denote the corresponding 
Novikov completion of the ring $\gL=\ZZZ G$.
The next theorem 
relates the homotopy type of the completed chain 
complex 
\bq\lb{f:univ_nov}
\wh\CC_*(\wi M,[\o])=
\wh\gL_{[\o]}
\tens{\gL} \CC_*(\wi M)
\end{equation}
(where $\CC_*(\wi M)$ denote the cellular chain complex
of the universal covering of $M$)
and the geometrical properties of the form
$\o$.
Denote by  $S_k(\o)$  the set of zeros of 
$\o$ of index $k$, and by $S(\o)$  the set of 
all zeros.
\beth\mlb{t:nov_complex}
There is a chain complex 
$\NN_*$
over the ring 
$\wh\gL_{[\o]}$
such that:
\beli\item 
$\NN_k$ 
is freely generated over $\wh\gL_{[\o]}$ by 
$S_k(\o)$.
\item 
The chain complexes 
$\NN_*$
and 
$~\wh\CC_*(M,[\o])$
are chain homotopy equivalent.
\enli
\enth
For the case of 
integral  classes $[\o]\in H^1(M,\ZZZ)$
the theorem follows from 
the existence of the Novikov complex
$\NN_*(f,v)$, associated to a circle-valued
Morse map $f$ such that $df=\o$
and any transverse $f$-gradient $v$
(see \cite{patou}). The case of rational classes
$[\o]\in H^1(M,\QQQ)$ follows immediately.
For the case of  1-forms 
belonging to arbitrary cohomology classes
see the later paper by F.Latour \cite{latour},
and also the works of  D.Sch\"utz
\cite{schuetz1}, \cite{schuetz2}.

\subsection{Twisted Novikov homology}
\mlb{su:twi_hom}

In this subsection 
we begin our study of the 
twisted Novikov homology.
Let $M$ be a connected  CW complex.
Let $\r:G\to GL_R(V)$
be any right representation,
$\pi:G\to H$  an epimorphism,
$\xi:H\to\RRR$ a \ho.
In the introduction 
we have associated to these data 
a chain complex 
$\wh \CC_*(M,\r_\pi,\xi)$
over $\Lxi$
(where $\L=R[H]$).
When  $\xi:H\to\RRR$
is a monomorphism and $R$ is a 
principal ideal domain,
the ring $\Lxi$
is also a principal ideal domain.
In this case put
\begin{gather*}\lb{f:def_nov}
\wh b_i(M,\r_\pi,\xi)
=
\rk_{\Lxi} H_i(\wh \CC_*(M;\r_\pi,\xi)),\\
\wh q_i(M,\r_\pi,\xi)
=
\tn_{\Lxi} H_i(\wh \CC_*(M;\r_\pi,\xi)).
\end{gather*}
(where $\tn$ stands for the torsion number 
of the module,
that is, the minimal possible number 
of generators 
of its torsion submodule).
Observe that if $R$ is a field, then all
the numbers $\wh q_i(M,\r_\pi,\xi)$
vanish.
\beth\mlb{t:siko2}
For a given right representation $\r$
and a given homomorphism $\pi:G\to H$
the numbers $b_k(M, \r_\pi,\xi)$
do not depend on the monomorphism $\xi$.
There is a set $\G\sbs H_\RRR'=Hom(H,\RRR)$
which is a 
finite union  of integral hyperplanes
such that in every connected component of 
the complement
$H_\RRR' \sm \G$ the numbers 
$\wh q_k(M,\r_\pi,\xi)$ do not depend on $\xi$.
\enth
In the particular case of 
the trivial representation
this theorem is due to J.-Cl.Sikorav 
(see \cite{pasbor});
the proof in the general 
case is similar.

\beco\mlb{c:vanishing_cones}
There is an open polyhedral conical subset 
$S\sbs H'_\RRR=\Hom(H,\RRR)$
such that 
a monomorphism $\xi:H\to \RRR$
is in $S$ if and only if 
the Novikov homology
$\wh H_*(M; \r_\pi,\xi)$ 
vanishes.
\enco

It is natural to ask, whether we can drop 
the condition of injectivity  of $\xi$
in the Corollary above 
and still keep the conclusion
of the Corollary.
I do not know if the answer 
is positive in general, 
but it is the case when $M$ 
is a compact  connected
3-manifold with $\chi(M)=0$
and $\pi:G\to H$ is  
the projection onto
$H=H_1(M,\ZZZ)/Tors$
(see Section \mrf{s:three}).

Let us proceed to the applications of 
the twisted Novikov homology 
in the theory of Morse forms.
Now we shall assume that $M$ is a 
closed connected $\smo$ manifold.
Let $\o$ be a Morse form on $M$.
The de Rham cohomology 
class $[\o]\in H^1(M,\RRR)$ of $\o$
can  be identified with a homomorphism
$G=\pi_1(M)\to\RRR$.
Let $\pi:G\to H$
be an epimorphism, such that $[\o]$ factors through $\pi$,
so that $[\o]=\xi\circ\pi$,
where $\xi:H\to\RRR$ is a homomorphism.

\beth\mlb{t:nov_complex_twi}
There is a chain complex $\NN^\rho_*$
over the ring $\Lxi$
such that:
\beli\item 
$\NN^\rho_k$ is a free $\Lxi$-module 
with $n\cdot m_k(\o)$ free generators.
\item 
The $\Lxi$-modules 
$H_*(\NN^\rho_*)$
and 
$\wh H_*(M; \r_\pi,\xi)$
are isomorphic.
\enli
\enth
\Prf
Let $\wh V^H=\Lxi\tens{\L} V^H$;
the composition of the right 
representation $\r_\pi$
with the natural inclusion
$GL_\L(V^H)\rInto GL_{\Lxi}(\wh V^H)$
determines a right representation 
$\wh\r_\pi$ of $G$ in $GL_{\Lxi}(\wh V^H)$.
Using the factorization $[\o]=\xi\circ\pi$
it is not difficult to check that
the  representation $\wh\r_\pi$ extends  
to a structure of a right
$\wh \gL_{[\o]}$-module on $\wh V^H$, and we have: 
\bq\lb{f:two_twisted}
\wh\CC_*(M; \r_\pi,\xi) =
\wh V^H\tens{\wh \gL_{[\o]} }\wh \CC_*(\wi M,[\o]).
\end{equation}
Now our theorem follows immediately from 
Theorem \mrf{t:nov_complex}. 
$\qs$

\beco\mlb{c:vanish}
Let $\o$ be a closed 1-form without zeros.
Assume that 
$$
[\o]=\xi\circ\pi
\in \Hom(G,\RRR) =
H^1(M,\RRR)
$$
where $\xi:H\to\RRR$ is a homomorphism,
and $\pi:G\to H$
is an epimorphism.
Then 
$\wh H_*(M; \r_\pi,\xi)=0.$
$\qs$
\enco

\beco\mlb{c:vanishh}
Let $\o$ be a Morse form.
Assume that 
$$
[\o]=\xi\circ\pi
\in \Hom(G,\RRR) =
H^1(M,\RRR)
$$
where $\xi:H\to\RRR$ is a monomorphism,
and $\pi:G\to H$
is an epimorphism.
Assume that $R$ is a principal ideal domain.
Then 
\bq
m_i(\o)
\geq \wh b_i(M,\r_\pi,\xi)
+
\wh q_i(M,\r_\pi,\xi)
+
\wh q_{i-1}(M,\r_\pi,\xi).
\hspace{1cm}\qs
\end{equation}
\end{coro}


\section{Fitting invariants of 
chain complexes}
\mlb{s:fitt_invar}

Let $A$ be a 
\fg~module over a commutative ring
$Q$. Let 
$$0\lTo A\lTo C_0\lTo^\l C_1$$
be a presentation for $A$, with $C_0,~C_1$
free \fg~modules.
By definition  the $k$-th Fitting ideal 
of $A$ 
is  
 the ideal generated by all
$(n-k)\times (n-k)$-minors of the matrix of the 
\ho~$\l$, where $n=\rk C_0$
(see \cite{eisenbud}, \S 20.2). 

In this section we give a generalization 
of this construction
and for every chain complex
$C_*$ of free \fg~ $Q$-modules
we define a family
of ideals of $Q$
({\it the Fitting ideals of $C_*$})
which are chain homotopy invariants 
of $C_*$. If $Q$ is a factorial ring,
each ideal has its greatest common divisor,
thus we derive from the family of 
the Fitting ideals 
a family of elements of $Q$,
which are called 
{\it the Fitting invariants of $C_*$}.
When $Q$ is a principal ideal domain, 
these invariants determine
the homology of $C_*$
(see Subsection \mrf{su:homol}).

If $C_*$ is a chain complex over a 
{\it non-commutative} ring
$\gL$, and $V$ is a left $Q$-module,
which is also a right $\gL$-module,
we can form the tensor product $V\tens{\gL} C_*$
and consider the Fitting invariants 
of the resultant complex. A particular 
case of this construction
leads to the well known {\it knot polynomials }
$\D_k(t)$ of \cite{crfox}, Ch. 8. We discuss this and 
similar constructions in 
Subsections \mrf{su:twi_mod} 
and \mrf{su:fitt_knots}.

While the definition 
of the Fitting ideals for chain complexes
is apparently new, 
many similar constructions 
exist already in the literature.
Let us mention for example the 
invariants of knots, deduced from
the representation spaces of 
the fundamental group of the knot
(see the paper \cite{Thang_repre} 
of Le Ty Quok Thang), and the twisted 
Alexander-Fox 
polynomials of V.Turaev
(see \cite{turaev_two_comple}).

\subsection{Matrices of homomorphisms: terminology}
\mlb{su:matrices}

This subsection is purely terminological:
we describe the conventions with which 
we shall be working.
Let $R$ be a ring (non-commutative in general).
Let $A,B$ be free finitely generated 
left modules over $R$.
Choose  a finite basis 
$\{e_i\}_{1\leq i\leq k}$ 
in $A$, and a finite basis 
$\{f_j\}_{1\leq j\leq  m}$ 
in $B$.
Let $\phi:A\to B$ be a module  homomorphism.
Write 
\bq\lb{f:matrix}
\phi(e_i)= \sum_j M_{ij} f_j.
\mxx{ with } M_{ij}\in R.
\end{equation}
The matrix $(M_{ij})$
will be denoted $ M(\phi)$
and called {\it the matrix of the \ho~$\phi$}
\wrt~the chosen bases.
Thus the coordinates of the images $\phi(e_i)$
of the basis elements of $A$
in $B$ are {\it the rows } of the matrix $M(\phi)$
(which  has $k$ rows and $m$ columns).
Here is the composition formula:
$$
M(\phi\circ\psi)
=
M(\psi)\cdot M(\phi),
$$
(where~ 
$\psi:C\to A,~\phi:A\to B$ are homomorphisms
of left modules, and $\cdot$ ~~stands for 
the usual matrix product).
This way of associating a matrix to
a module homomorphism will be called 
{\it row-wise}. For a free module 
$A$ with $k$ free generators
the map $$
\phi\mapsto M(\phi): \Hom(A,A)\to Mat(k\times k, R)
$$
is therefore an anti-homomorphism. 

In many recent textbooks on linear algebra one 
finds  another convention:
\bq\lb{f:matrix_column}
\phi(e_i)= \sum_j \wi M_{ji} f_i,
\end{equation}
so that $\wi M(\phi) = M(\phi)^T$, and
the coordinates of the images of the 
basis elements of $A$
in $B$ form  the columns of the matrix $\wi M$.
This way of associating the matrix to
a module homomorphism will be called 
{\it column-wise}.

When the ring $R$ is commutative the convention 
\rrf{f:matrix_column}
leads to the following composition formula:
$$
\wi M(\phi\circ\psi)
=
\wi M(\phi)\cdot \wi M(\psi),
$$
therefore  the map 
$$\phi\mapsto \wi M(\phi),\quad Hom(A, A)\to Mat(k\times k , R)$$
is a ring homomorphism.

\subsection{The Fitting invariants of 
chain complexes over 
commutative rings}
\mlb{su:determ_commut}

\bede\mlb{d:regular}
A chain complex 
$$
C_*
=\{\dots \lTo C_k \lTo^{\pr_{k+1}} C_{k+1}\dots \}
$$
of left modules over some ring $Q$
is called {\it regular} if 
every $C_i$ is a 
\fg~free $Q$-module, and
$C_i=0$ for $i<0$.
\end{defi}
In this subsection 
$Q$ is a commutative factorial  ring.

\bede
Let $C_*$ be a regular 
chain complex of $Q$-modules, and 
$k\in\NNN$.
Choose any finite bases in the modules
 $C_k, C_{k+1}$
and let $M(\pr_{k+1})$ be the matrix of 
$\pr_{k+1}$ 
\wrt~these bases.

 Let
$I_s(\pr_{k+1})$
denote 
 the ideal in $Q$ generated by all $s\times s$-minors 
of $M(\pr_{k+1})$.
(Here we assume that $s$ is an integer with 
$0<s\leq \min( \rk C_k, \rk C_{k+1})$. 
 If $s > \min( \rk C_k, \rk C_{k+1})$, 
then put by definition  
$I_s(\pr_{k+1})=0$, and for $s\leq 0$
 put $I_s(\pr_{k+1})=Q$.)
\end{defi}

The next lemma is a well-known consequence 
of the Binet-Cauchy formula (see for example
 \cite{Mcd}, p. 25).
\bele\mlb{l:det_indep}
The ideal $I_s(\pr_{k+1})$ does not 
depend on the particular choice of
bases in $C_k$ and $C_{k+1}.\qs$
\enle

It is clear that the ideal $I_s(\pr_{k+1})$
is not in general an invariant 
of the homotopy type of
the chain complex. However we can 
re-index the sequence $I_s(\pr_{k+1})$
and obtain homotopy invariants.

\bede\mlb{d:def_invar}
Put
$$
J^{(k)}_m(C_*)
=
I_{\rk C_k - \rk C_{k-1} - m+1} (\pr_{k+1}).
\quad\footnote{
This re-indexing may seem 
arbitrary, but we shall 
see that it fits
with the usual notation 
for the knot polynomials.
})
$$
\end{defi}
\bepr\mlb{p:det_invar}
For every $m, k$ the ideal 
$J^{(k)}_m(C_*)$
is a homotopy invariant of the regular 
chain complex $C_*$.
\enpr
\Prf
Let $F$ be a free \fg~$Q$-module.
Let $T_*(i,F)$ denote the chain complex
$$
 T_*(i,F) 
=
\{0 {\lTo} ~{\cdots}~ 0 \lTo  
F \lTo^{\id} F {\lTo} 0 {\lTo}  {\cdots}\}
$$
concentrated in degrees $i, i+1$.
A chain complex isomorphic to a direct sum 
of complexes $T_*(i,F)$ for some $i\geq 0$
and some $F$ will be called {\it trivial}.
The next lemma is one of the versions of the
Cockroft-Swan theorem \cite{CockroftSwan},
the proof is similar to the
the proof of Lemma 1.8 in \cite{pasur}.
\bele\mlb{l:thickenings}
Let $C_*, D_*$
be  chain homotopy equivalent  complexes.
Then there are  trivial 
chain complexes $T_*, T'_*$
such that
$C_*\oplus T_*\approx D_*\oplus T'_*$.
\enle
\Prf 
Let 
$$K_*=\bigoplus_{i\in\ZZZ}
T(i,C_i).
$$
This is a  trivial chain complex. Let 
$\phi:C_*\to D_*$ be 
a chain homotopy equivalence.
Define a chain map $\psi: C_*\to D_*\oplus K_*$
as follows:
\begin{gather}
\psi(c)= (\phi(c),c, \pr c) \in D_i\oplus 
C_i\oplus C_{i-1} 
\mxx{ for } c\in C_i.
\end{gather}
It is clear that $\psi$ is a chain 
homotopy equivalence 
which is a split monomorphism.
The quotient chain complex
$(D_*\oplus K_*)/\Im \psi$
is a contractible chain complex of free \fg~
modules.
Thus we obtain
an exact sequence
\bq\lb{f:split_exact}
0 \rTo C_*\rTo^\psi D_* \oplus K_* \rTo S_* \rTo 0,
\end{equation}
where $S_*$ is a regular  acyclic chain complex.
Such a sequence splits (see for 
example \cite{cohen}, 13.2),
and we obtain an isomorphism 
$$
C_*\oplus S_*
\approx 
D_*\oplus K_*.
$$
It is easy to prove that 
there is a free  trivial chain complex
$R_*$, such that $S_*\oplus R_*$ is trivial, and 
this completes the proof of the lemma
with $T_*=S_*\oplus R_*,~ T'_*=K_*\oplus R_*$.
 $\qs$

Now let us return to the 
proof of our proposition.
In view of the preceding lemma it 
suffices to check the 
following easily proved 
assertion:
the ideal 
$J^{(k)}_m(C_*)$
does not change if we add to the
$C_*$ the  chain complex
$T_*(i, F)$ where $F$ is a \fg~free
$Q$-module, and $i$ equals one
of the numbers $k-1, k$ or $k+1$.
$\qs$

\bede\mlb{d:det_id_det_invar}
The ideal 
$J^{(k)}_m(C_*)$ is  called 
{\it the $m$-th Fitting ideal of $C_*$}.
The GCD ({\it = the greatest common divisor})
of all the non-zero elements in the ideal
$J^{(k)}_m(C_*)$ will be denoted by $F^{(k)}_m(C_*)$
and called {\it Fitting 
invariant of $C_*$}.
This element 
is well-defined up to multiplication 
by invertible elements of $Q$.
The sequence 
\bq\lb{f:seq_fit}
\dots J_m^{(k)}(C_*) \sbs J_{m+1}^{(k)}(C_*)\dots 
\end{equation}
of the Fitting ideals
will be called {\it the Fitting 
sequence of $C_*$ in degree $k$}.
The subsequence 
of  \rrf{f:seq_fit}
formed by all non-trivial ideals,
is called
{\it reduced Fitting sequence 
of $C_*$ in degree $k$.}
\end{defi} 
Recall that an ideal $I\sbs Q$
is called {\it non-trivial}
if $I\not=0,~I\not= Q$.

Observe that the length of the reduced 
Fitting sequence in degree $k$
is not more than $\rk C_k$. 
\bele\mlb{l:localiz}
Let  $S\sbs Q^\bu$
be a multiplicative subset of $Q$.
Then up to invertible elements of $ S^{-1}Q $
we have:
$$
\fmkk {C_*}
=
\fmkk { S^{-1}C_*}.\hspace{1cm} 
$$
\enle
\Prf
It suffices to recall that
the GCD of elements of a factorial ring 
does not change when the ring is localized.
$\qs$

\bere\mlb{r:two_reg}
The Fitting invariants  $F_m^{(k)}$
for $k\leq 2$ can be defined in 
a slightly more 
general framework: 
\bede\mlb{d:reg}
A chain complex $C_*$
of left $Q$-modules 
will be called {\it 2-regular}
if $C_i=0$ for $i > 0$ and 
$C_0, C_1, C_2$ are \fg~
free $Q$-modules.
\end{defi}
Using the same procedure as above
we can define the Fitting invariants 
$F_m^{(i)}$
with  $i\leq 2$ for any 2-regular 
chain complex over
a commutative ring $Q$. 
\enre

\subsection{The case when $Q$ is principal}
\mlb{su:homol}

Assume that 
$Q$ is a principal ideal domain.
We shall show that 
in this case the sequence of Fitting invariants
determines the homology 
of the chain complex.
We shall give only the statements 
of the theorems; 
the proofs are obtained
applying the standard results about 
the structure of the modules 
over principal ideal domains.

Let us begin with a theorem which shows how to
compute the Betti
numbers $$
b_k(C_*)
=
\rk_Q H_k(C_*)
$$
and the torsion numbers
$$
q_k(C_*)=
\tn_Q H_k(C_*) 
$$
from the Fitting invariants.
(Recall that the {\it torsion number of a module $X$}
is the minimal possible number of 
generators of the torsion submodule of $X$.)
Let $C_*$ be a regular chain complex
of  $Q$-modules, let
$\g_k=\rk C_k$.
Consider the subsequence of the 
Fitting sequence starting with 
$J^{(k)}_{-\g_{k-1}+1}=I_{\g_k }(\pr_{k+1})$:
$$
J_{-\g_{k-1}+1}, ~J_{-\g_{k-1}+2}, \dots 
$$
Let $A_k$ be the number of zero 
ideals in this sequence.
Let $\varkappa_k$ be the cardinality of 
the reduced Fitting sequence
of $C_*$ in degree $k$.
\beth\mlb{t:numbers}
We have
\begin{align}\lb{f:numbers}
  b_k(C_*)  = & A_k+A_{k-1} -\g_{k-1}; & &
\\
  q_k(C_*)  = & \varkappa_k.  & &\qs
\end{align}
\enth
The torsion submodule 
of the homology is also determined
by the Fitting invariants.
Write the reduced Fitting 
sequence in degree $k$ as follows:
\bq\lb{f:reduced_fit}
I_1,\dots, I_{\varkappa_k}.
\end{equation}
Let 
$\t_s$ be the GCD of all the elements of $I_s$
then 
we have:  $\t_{s+1}~|~\t_{s}$ for every $s$.
Put $\l_s= \t_{s}/\t_{s+1}$.

\beth\mlb{t:torsion}
The elements $\l_s\in Q$ are non-invertible,
for every $s$ we have:
$\l_{s+1}~|~\l_{s}$
and 
$$
Tors~ H_k(C_*)
\approx
\bigoplus_{i=1}^{\varkappa_k} Q/\l_i Q.\hspace{1cm}\qs
$$
\enth

\subsection{The twisted Fitting  invariants 
of $\ZZZ G$-complexes}
\mlb{su:twi_mod}

In this subsection we apply the 
Fitting invariants of the preceding subsection
to construct invariants
of chain complexes over non-commutative rings.
The most important for us are modules over
group rings, and we limit ourselves to this case, 
although there are obvious generalizations.

Let $G$ be a group and  let
$C_*$ be a regular chain complex 
over $\gL = \ZZZ G$. Let $\t:G\to GL_Q(W)$
be a right representation,
where $W$ is a \fg~free $Q$-module
over some commutative factorial ring $Q$.
Form the tensor product
$$
C_*(\t)
=
W\tens{Q} C_*;
$$
the Fitting invariants
$\fmkk{W\tens{Q} C_*}$ of this complex
will be denoted by
$\d_m^{(k)}(C_*, \t)$ and called 
{\it twisted Fitting invariants of $C_*$}.

\bere
\mlb{r:two_regg}
Similarly to Remark \mrf{r:two_reg}
we obtain the Fitting invariants 
$\d_m^{(k)}(C_*,\t)\in Q$ where $k\leq 2$ and 
$C_*$ is a 2-regular chain complex.
\enre

In this paper we shall be interested  
mainly in the case when 
the ring $Q$ is the group ring 
of some free abelian group,
or a Novikov completion 
of such group ring, or a 
localization of such group ring.

Let $H$ be a free abelian \fg~
group, and $R$ a commutative factorial
ring. The group ring $R[H]$ will be denoted by 
$\L$. Let $V$
be  a \fg~free left $R$-module.
Let $\r:G\to GL_R(V)$
be a right representation of $G$
and $\pi:G\to H$ is an 
epimorphism.
Recall that in Introduction
we have associated to this data 
a right representation
$\r_\pi:G\to GL(V^H)$
where $V^H=\L\tens{R} V$.

\bede\mlb{d:twi_fitt}
Let $C_*$ be a regular chain complex
over $\ZZZ G$.
The Fitting invariant
$\d_m^{(k)}(C_*,\r_\pi)\in \L$
will be 
called {\it the twisted 
Fitting invariant of
$C_*$ \wrt~$(\r,\pi)$.}
\end{defi}

In the rest of this subsection
we investigate the behaviour  
of the Fitting
invariants of $C_*$
\wrt~ certain completions and 
localizations
of the representation 
$\r_\pi$.

\bede\mlb{d:monicc}
Let $\xi:H\to\RRR$
be a non-zero \ho. 
An element $x\in \L=R[H]$
is called {\it $\xi$-monic}
if 
\begin{gather*}
x=x_0 h_0+\sum_{i=1}^s x_ih_i
\mxx{ with } x_i\in R, h_i\in H, \\
\mxx{ where } x_0\in R^\bu
\mxx{ and }
\xi(h_i) <\xi(h_0) 
\mxx{ for every } i\not=0.
\end{gather*}
The multiplicative subset of 
all $\xi$-monic 
elements will be denoted $S_\xi$.
The ring 
$S_\xi^{-1}\L$ will be also denoted by
$\Lx$.
\end{defi}
The next proposition is immediate.
\bepr\mlb{p:monic_inver}
An element $x\in \L$
is $\xi$-monic if and only if it is 
invertible in $\Lxi$. $\qs$
\enpr
Therefore the ring $\Lx$ can be 
considered as a 
subring of $\Lxi$.
Let $\wi\r_{\pi,\xi}$
denote the composition 
$$
G\rTo{\r_\pi} GL_\L(V^H){\rInto} 
GL_{\Lx} (\Lx\tens{\L}V^H).
$$
This is a right representation 
of $G$, and we obtain the 
corresponding twisted Fitting 
invariants 
$\d_m^{(k)} (C_*,\wi\r_{\pi,\xi}) $.
It is clear that
\bq\lb{f:fitt_local}
\d_m^{(k)} (C_*,\r_{\pi,\xi})
=
\d_m^{(k)} (C_*,\wi\r_{\pi,\xi}). 
\end{equation}
Similarly, let $\wh\r_{\pi,\xi}$
denote the composition 
$$
G\rTo{\r_\pi} GL_\L(V^H){\rInto} 
GL_{\Lxi} (\Lxi\tens{\L}V^H).
$$
This is a right representation 
of $G$, and we obtain the 
corresponding twisted Fitting 
invariants 
$\d_m^{(k)} (C_*,\wh\r_{\pi,\xi}) $.

In the rest of this subsection we restrict
ourselves to the particular case $R=\ZZZ$,
although some of the results 
can be proved in a more general setting.

\bepr\mlb{p:twi_fitt_nov}
Let $\xi:H\to\RRR$
be a monomorphism. Let $R=\ZZZ$. 
Then 
$$
\d_m^{(k)}(C_*,\r_\pi)
=
\d_m^{(k)}( C_*,\wh\r_{\pi,\xi} ). 
$$
\enpr
\Prf 
We shall reduce the proof to the 
equality \rrf{f:fitt_local}.
\bepr[\cite{pasbor}]\mlb{p:mono_principal}
If $\xi:H\to\RRR$
is monomorphic,
then $\Lx$ is a principal ideal domain.
\enpr

\beco\mlb{c:gcd}
For every two elements $a,b\in \L$
we have
$$
GCD_\L(a,b)
=
GCD_{\Lx}(a,b)
=
GCD_{\Lxi}(a,b).\hspace{1cm} \qs
$$
\enco
The proposition 
follows immediately. $\qs$

Now we can explain how to compute
the twisted Novikov homology 
of a chain complex in terms of its 
twisted Fitting invariants.
Let  $C_*$ be a regular 
chain complex over $\gL=\ZZZ G$.
Let $\pi:G\to H$ an epimorphism onto 
a free \fg~abelian group, and 
$\r:G\to GL(V)$ be a right representation,
where $V\approx \ZZZ^n$
is a free \fg~module over $\ZZZ$. 
Let $\xi:H\to\RRR$ be a monomorphism,
so that $\Lxi$ is a principal ideal domain.
Applying the results of Subsection 
\mrf{su:homol}
we obtain the following description
of the twisted Novikov homology. 

Let $\g_k=n\cdot\rk C_k$, 
where $n=\rk_\ZZZ V$, and $\rk C_k$
is the number of free generators of 
$C_k$. 
Consider the segment of
the Fitting sequence 
of $C_*(\r_\pi)=V^H\tens{\L} C_* $
in degree $k$
starting with 
$J^{(k)}_{-\g_{k-1}+1}=I_{\g_k }(\pr_{k+1})$:
$$
J_{-\g_{k-1}+1}, ~J_{-\g_{k-1}+2},\  \dots 
$$
Let $A_k$ be the number of zero 
ideals in this segment.
Let $I_1\sbs \dots \sbs I_{B_k}$
be the reduced Fitting sequence
of $C_*(\r_\pi)$
in degree $k$,
let $\l_s\in \L$
be the GCD of the ideal $I_s$.
Then $\l_i~|~\l_j$ for $i\geq j$.
Let $\varkappa_k(\xi)$
be the number of {\it non-$\xi$-monic}
elements $\l_j$.

\beth\mlb{t:numberss}
We have
\bega\lb{f:numberss}
\wh b_k(X,\r_\pi,\xi) 
= 
A_k+A_{k-1} -\g_{k-1};
\\
\wh q_k(X,\r_\pi,\xi)
=
\varkappa_k(\xi).\hspace{1cm}\qs
\end{gather}
\enth
As for the 
torsion submodule in degree $k$, 
let $\t_s=\l_s/\l_{s+1}\in \L$.
\beth\mlb{t:torsionn}
The elements $\t_s\in \Lx$ are 
non-invertible.  For every $s$ 
$$
\t_{s+1}\mxx{ divides } \t_{s}\mxx{ in } \Lx
$$
and 
$$
Tors~ \wh H_k(X,\r_\pi,\xi)
\approx
\bigoplus_{i=1}^{\varkappa_k(\xi) } \Lxi/\t_i \Lxi.
\hspace{1cm}\qs
$$
\enth

\subsection{The twisted Fitting
invariants 
of $\ZZZ G$-modules}
\mlb{su:twi_modd}

Let $G$ be a group, and 
$N$ be a left $\ZZZ G$-module, 
admitting a free \fg~  resolution
$$
0\lTo R_0\lTo R_1\lTo R_2 {\lTo} ...
$$
over $\ZZZ G$
(so that the homology of $R_*$ 
vanishes in all dimensions except 
zero, and $H_0(R_*)\approx N$.)
Let  $Q$ be a commutative factorial ring, 
and $V$ a finitely generated 
free left $Q$-module.
Let $\t:G\to GL(V)$ 
be a right representation 
of $G$.
Consider the  chain complex 
\bq\lb{f:tensorr}
V\tens{\ZZZ G} R_*
\end{equation}
of free left $Q$-modules
(the module $V$ is endowed 
with the structure of a right 
$\ZZZ G$-module 
via the representation
$\t$).
Observe that the twisted Fitting invariants
\bq\lb{f:det_invarr}
\fmkk{V\tens{\ZZZ G} R_*}
\end{equation}
depend only on $N$ and $\t$,
but not on the particular choice 
of the resolution $R_*$
(indeed, any two resolutions are 
chain homotopy equivalent).
The most important for us are 
the Fitting invariants
corresponding to the second 
boundary operator.
\bede\mlb{d:det_groups}
The  element
$
F_m^{(2)} (V\tens{\ZZZ G} R_*)
\in Q
$
will be called {\it $m$-th Fitting  invariant}
(or {\it $det$-invariant}) of the pair $(G,\t)$
and denoted by
$$
\d_m(G,N,\t)
=
F_m^{(2)} (V\tens{\ZZZ G} R_*)
\in Q.
$$
This element is well-defined 
up to multiplication by an invertible 
element of the ring $Q$.
\end{defi}

Similarly,
we obtain the Fitting invariants 
$\d_m(G,N,\t)$
for the case when $N$ is a free left $\ZZZ G$-module
which admits a 2-regular resolution.
(see Remarks \mrf{r:two_reg},
 ~\mrf{r:two_regg}).

\begin{exam}\mlb{e:knot_poly}
Let $K$ be an oriented  knot; put 
$G=\pi_1(S^3\sm K)$.
Let $N$ be an open tubular \nei~of $K$.
Choose any finite CW decomposition 
of $S^3\sm N$.
The cellular chain complex
$\CC_*(\wi{S^3\sm N} )$
is a free $\ZZZ G$-resolution of 
the module 
$H_0(\wi{S^3\sm N})\approx \ZZZ$.
The canonical epimorphism $\ve:G\to\ZZZ$, 
sending each positively oriented
meridian to $1$
extends to a ring  
\ho~$\t:\ZZZ G\to(\ZZZ[\ZZZ])^\bu=GL(1,\ZZZ[\ZZZ] ) $.
The Fitting invariant $\d_1(G,\ZZZ, \t)\in \ZZZ[\ZZZ]$
equals clearly the Alexander
polynomial of the knot.
More generally,
$\d_i(G,\ZZZ, \t)$
is the knot polynomial $\D_i(t)$
(in the terminology of \cite{crfox}, Ch. 8).
\end{exam}

\mlb{e:mod_poly}
The previous example has a 
natural generalisation.
Let $G$ be a finitely presented group.
The abelian group $\ZZZ$ endowed with 
the trivial action of $G$
admits a 2-regular free resolution. 
Let $V$ be a free \fg~left $R$-module,
where $R$ is a commutative factorial ring, 
and $\r:G\to GL_R(V)$ be a right 
representation.
Let $\pi:G\to H$ be a homomorphism 
of $G$ to a free abelian finitely 
generated group. 
We have then the twisted 
Fitting invariants 
corresponding to the right 
representation
$\r_\pi:G\to GL_\L(V^H)$,
where $\L=R[H]$.
\bede\mlb{d:fitt_modules}
The twisted Fitting invariant 
$\d_m(G,\ZZZ,\r_\pi)\in \L$
will be denoted also by 
$\d_m(G,\r_\pi)$.
The first Fitting invariant
$\d_1(G,\ZZZ,\r_\pi)$
will be   denoted  also by $A(G, \rho_\pi)$.
\end{defi}
It turns out that 
the  Fitting invariants with 
non-positive indices vanish:
\bele\mlb{l:vanish_dets}
$$
\d_m(G,\r_\pi)=0
\mxx{ for } m\leq 0.
$$
\enle
\Prf
Pick a presentation of a group $G$, 
let $g_1, ..., g_s$ be the generators,
and $r_1, ..., r_l$ be the relations.
Write a 2-regular free resolution 
for $\ZZZ$
over $\gL=\ZZZ G$ as follows:
$$
F_*=\{ 0  \lTo \gL \lTo^{\pr_1} \gL^s
\lTo^{\pr_2}   \gL^l{\lTo } ...\};
$$
here the free generators $e_1,...,e_s$
of the module $F_1=\gL^s$ correspond 
to the generators $g_1,...,g_s$
of $G$ and the homomorphism
$\pr_1$ is given by 
$\pr_1(e_i)=1-g_i$.
We can  assume that $l\geq s$, and that 
$t=\pi(g_1)$ is one of the
 free generators of the abelian group $H$.
Let us consider  the image of 
our  Fitting invariant
in the field of fractions $\RR$ of $R[H]$.
This element coincides with the Fitting invariant 
of the chain complex 
$$\FF_* =
\RR^n\tens{\ZZZ G} F_*
=\{ 0  \lTo \RR^n\lTo^{\wi\pr_1} \RR^{ns}
\lTo^{\wi\pr_2}
\RR^{nl}{\lTo} ...\}.
$$
Observe that $\wi\pr_1$ is an 
epimorphism. Indeed, the restriction 
of $\wi\pr_1$ to the first 
direct summand $\RR^n$ of $\RR^{ns}$
equals
$1-t\r(g_1):\RR^n\to\RR^n$, and the determinant 
of this map is non-zero, 
therefore invertible in $\RR$.
Thus the rank of the matrix $\wi\pr_2$
is not more than 
$n(s-1)$, and the first Fitting invariant
which can be non-zero
is
$$
J_1^{(2)}(\FF_*)
=
I_{\rk \FF_1 -\rk \FF_0 }(\pr_2).
$$
 $\qs$

The Fitting invariants of modules
are useful for computations of 
the homology with twisted coefficients.
Let $X$ be a connected finite CW complex,
put $G=\pi_1(X)$.  
Let $\pi:G\to H$ an epimorphism onto 
a free \fg~abelian group,
$\r:G\to GL(V)$ be a right representation,
where $V\approx \ZZZ^n$
is a free \fg~module over $\ZZZ$.
Let $\CC_*(\wi X)$ denote the cellular chain
complex of the universal covering of $X$. 

\bepr\mlb{p:groups_CW}
$\d_i(\CC_*(\wi X),\r_\pi) =\d_i(G,\r_\pi)$.
\enpr
\Prf
We have:
$$
H_0(\wi X)\approx \ZZZ, \quad H_1(\wi X) =0.
$$
The module $H_2(\wi X)$ can be non-zero, 
and therefore
the chain complex $\CC_*(\wi X)$
fails in general to be a resolution of 
the module $\ZZZ$.
Choose any subset $S\sbs \CC_2(\wi X)$
generating the $\ZZZ G$-submodule $Z_2(\wi X)$
of 2-cycles in the complex $\CC_*(\wi X)$.
Put $\gL=\ZZZ G$ and 
let $\gL^S$ be the free $\gL$-module 
generated by the set $S$.
Extend the identity map $S\rTo^{\id } S$
to a $\gL$-module map
$\phi:\gL^S\to \CC_2(\wi X)$.
Put 
$$\CC_3'=\CC_3(\wi X)\oplus \gL^S,\quad
\pr'_3=(\pr_3,\phi):\CC'_3\to \CC_2(\wi X).
$$
The chain complex $\CC_*'$ is 2-regular and 
its second homology vanishes.
Applying the same procedure 
to the third, fourth, etc. homology modules,
we obtain an acyclic 2-regular resolution 
$D_*$ of the module $\ZZZ$
such that
$$
\CC_*(\wi X)\sbs D_* \mxx{ and }
D_i/\CC_i(\wi X)=0 
\mxx{ for } i\leq 2.
$$
Thus 
$$
\d_i(\CC_*(\wi X),\r_\pi)=
\d_i(D_*,\r_\pi);
\quad
\d_i(G,\r_\pi)=
\d_i(D_*,\r_\pi).
\hspace{1cm}
\qs$$

\subsection{Fitting invariants of knots and links}
\mlb{su:fitt_knots}

Let $L$ be an oriented link.
Put $G=\pi_1(S^3\sm L)$.
The group $G$ is finitely 
presented, therefore the module 
$\ZZZ$ has a 2-regular free resolution over $\ZZZ G$.
Thus  for any epimorphism
~$\pi:G\to H$
and any right representation 
$\r:G\to GL_R(V)$
we obtain a sequence of elements
$$
\d_1(G,\r_\pi),~ \d_2(G,\r_\pi),... \in \L=R[H],
$$
defined up to multiplication by an invertible element of $R$.
\bede\mlb{d:fitt_links}
The elements 
$$\d_i(G,\r_\pi)=\d_i(\pi_1(S^3\sm L),\r_\pi)$$
are called {\it the Fitting 
invariants of the link $L$}.
\end{defi}

We shall discuss these invariants 
and their relations with Novikov homology
in more details 
in Section \mrf{s:three}.


\section{Twisted Alexander polynomials}
\mlb{s:det_twialex}

We begin with  a recollection of 
M.Wada's definition of the 
twisted Alexander polynomial;
this occupies the first two subsections.
In the rest of the section we discuss  
the relations between
the  Fitting invariants 
and twisted Alexander polynomials.
The  one-variable case
and the  multi-variable case 
are slightly different from each other
and are considered separately.

\subsection{W-invariant of a matrix}
\mlb{su:w_invar}
Let $\gL$ be a ring with a
unit (non-commutative in general);
let $\l:\gL\to Mat(n\times n,\L)$
be a ring homomorphism, where $\L$ is 
a commutative factorial ring.

For any  matrix $\AA$ over 
$\gL$ we can substitute 
 instead of each of its
matrix entries $\AA_{ij}$
its image \wrt~
$\l$. 
The result of this operation 
will be denoted by $\psi(\AA)$;
the size of $\psi(\AA)$
is $n$ times the size of $\AA$.
We have:
$$
\psi(\AA_1\AA_2)=\psi(\AA_1)\psi(\AA_2),
$$
if the number of rows of $\AA_2$
equals the number of columns of $\AA_1$.

Now let $\BB$ be an $l\times s$ - matrix 
with coefficients
in $\gL$. Assume that $l\geq s-1$.
Let 
$$
\a=
\begin{pmatrix}
a_1 \\
\vdots \\
a_s \\
\end{pmatrix}
$$
be a column   of elements of $\gL$,
such that 
\bq\lb{f:product_zero}
\BB\cdot \a=0.
\end{equation}
We are going to associate
to these data an element
$$
W=W(\BB,\a,\l)\in \L.
$$

For an integer $j$ with $1\leq j\leq s$
denote by $\BB_j$ the $l\times (s-1)$-matrix
obtained from $\BB$ by suppressing the
$j$-th column. 
We have the $nl\times n(s-1)$-matrix $\psi(\BB_j)$.
Let $S=i_1<i_2< ... < i_{n(s-1)}$
be a sequence of integers in $[1, nl]$.
Let $\psi(\BB_j)^S$
denote the square matrix formed by 
all the matrix entries 
of $\psi(\BB_j)$ contained in  the rows 
with indices in $S$.
It is easy to deduce from
the condition
\rrf{f:product_zero}
that for every $i,j$ with $1\leq i,j\leq s$
we have
\bq\lb{f:equall}
\det\Big(\psi(\BB_j)^S\Big)
\cdot
\det\big(~\psi(a_i)~\big)
=
\det\Big(\psi(\BB_i)^S\Big)
\cdot
\det\big(~\psi(a_j)~\big).
\end{equation}
Let $Q_j(\BB)$ denote the GCD of 
the elements 
$\det\Big(\psi(\BB_j)^S\Big)$
over all $S$.
This element is defined up to 
multiplication by an invertible element
of $\L$.
\bede\mlb{d:w_invar}
Assume that
there exists $j$ with
$\det\big(~\psi(a_j)~\big)\not=0$.
The element 
$$
W(\BB,\a,\l)=\frac {Q_j(\BB)}{\det\big(~\psi(a_j)~\big)}
$$
of the fraction field of $\L$ 
will be called
W-invariant of the matrix $\BB$.

If $l\leq s-1$ we set by definition
$W(\BB,\a,\l)=0$.
\end{defi}

\subsection{Twisted Alexander 
polynomials: definition}
\mlb{su:twi_alex_def}

Now we apply the construction of 
the previous subsection
to define the \talp.
Let $G$ be a finitely presented group. Let 
$\pi:G\to H$ be an epimorphism of $G$ onto 
a finitely generated free abelian group $H$.
Let $R$ be a commutative factorial ring, and 
let $\l:G\to GL(n,R)$ be a group homomorphism.
The tensor product of this \ho~
with $\pi: G\to R[H]^\bu$
gives rise  a ring \ho~
$\l_\pi:\gL=\ZZZ G\to Mat(n\times n, \L)$,
where $\L=R[H]$.

Pick a finite presentation 
$p=(g_1,...,g_s~|~h_1,...h_l)$
of the group $G$,
with generators $g_j$
and relators $h_i$.
Let $\chape {h_i}{g_j}\in \ZZZ G$
denote the corresponding Fox derivative.
\bede\mlb{d:alex_matrix}
The matrix 
$$
\AA=\AA(G,p)=\Big(\chape {h_i}{g_j}\Big)
\in Mat(l\times s,\gL)
$$
will be called {\it Alexander matrix}
of the presentation 
$p$.
\end{defi}
Let 
$$
\a=
\begin{pmatrix}
1-g_1\\
\vdots\\
1-g_s\\
\end{pmatrix}
\in \gL^s.
$$
We have:
$$
\AA\cdot\a=0,
$$
therefore the constructions of the 
previous subsection
apply and we obtain the W-invariant
$W(\AA,\l_\pi,\a)$
in the fraction field of 
$R[H]$. The ring $R[H]$
being isomorphic to the ring of Laurent
polynomials in $k$ variables 
with coefficients in $R$ (where $k=\rk H$),
the W-invariant can be considered as 
a rational function of 
$k$ variables with coefficients 
in the fraction field
of $R[H]$.
M.Wada proves in \cite{Wada}
that this element does not depend on 
the particular choice
of the presentation $p$
(up to multiplication by an invertible element
of $\L$) and is therefore determined by
$\l$ and $\pi$. 
\bede\mlb{d:twi_al}
The element 
$W(\AA,\l_\pi,\a)$ of the fraction field of $R[H]$
is denoted by
$\D_{G,\l}$ and called {\it \talp~
of $G$ associated to $\l$ and $\pi$}.
\end{defi} 

\subsection{Relation with 
the Fitting invariants: 
the one-variable case}
\mlb{su:rel_det_one}

In the case when $H\approx \ZZZ$
the ring $\L=R[H]$ is isomorphic to $R[t, t^{-1}]$.
By definition 
the twisted Alexander polynomial is 
an element of the field of fractions of the ring $R[t]$.
Consider a multiplicative subset $\Sigma\sbs \L$ 
consisting  of all Laurent polynomials 
of the form
$a_it^i+ ... + a_jt^j$,
where $i,j\in\ZZZ,~i\leq j$ and $a_i, a_j$ are 
invertible elements of $R$.
\bepr\mlb{p:twial_local}
$\D_{G,\l}(t)\in \Sigma^{-1} \L$.
\enpr
\Prf 
We can assume that the first generator $g_1$
satisfies $\pi(g_1)=t\in R[t, t^{-1}]$.
The first coefficient of the polynomial 
$\det(1-t\l(g_1))$
equals  1 and the last 
coefficient equals $\pm \det(\rho(g_1))$,
which is an invertible element of $R$. $\qs$

Now we can proceed to the comparison of 
Fitting invariants and twisted Alexander polynomials.
Let $\r:G\to GL_R(V)$
be a right representation of $G$, where $V$ 
is a finitely generated free left $R$-module.
Choose a basis in $V$, then
we obtain a homomorphism
$$\ove{\r}:G\to GL(n,R);
\quad
\ove{\r}(g)=M(\r(g)).$$

\bepr\mlb{p:al_det_first}
The Fitting invariant 
$A(G,\rho_\pi)$
divides the \talp~$\D_{G,\ove{\r}}$
in the ring 
$\Sigma^{-1} \L$.
\enpr
\Prf
Let $g_1,...,g_s$ be generators of $G$ and
$h_1, ..., h_l$ be relators.
The free resolution of $\ZZZ$ over $\gL=\ZZZ G$
can be constructed as follows:
\bq\lb{f:pfff}
F_*=\{0\lTo \gL\lTo^{\pr_1}  
\gL^s\lTo^{\pr_2}  
\gL^l {\lTo}...
\end{equation}
where 
$$
M(\pr_1)
=
\begin{pmatrix}
1-g_1\\
\vdots\\
1-g_s\\
\end{pmatrix}
$$
and  $M(\pr_2)$
equals the Alexander matrix $\AA=\AA(G,p)$
corresponding to the presentation
$p=(g_1,\dots, g_s~|~h_1,\dots, h_l)$.
(see \cite{bz}, Ch. 9).
The polynomial $A(G,\r_\pi)$
is computed from the chain complex
$\FF_*=\L^n\tens{\ZZZ G} F_*$.
The matrix $M(\wi\pr_2)$ 
of the second boundary operator
in the chain complex $\FF_*$
is equal to the matrix $\psi(\AA)$
of Subsection \mrf{su:w_invar},
where $\psi=\ove{\r_\pi}$.
Now the proposition follows easily, 
since for every $j$ the 
element $Q_j(\AA)$ from Definition 
\mrf{d:w_invar}
is the GCD of a certain family
of $(s-1)n\times (s-1)n$-minors 
of $M(\wi\pr_2)$, and $A(G,\r_\pi)$
is the GCD of {\it all} 
$(s-1)n\times (s-1)n$-minors. $\qs$

\bepr\mlb{p:alex_det_links_one}
If the group $G$ has a presentation 
with $s$ generators and $s-1$ relations
(for example, $G$ is the fundamental group of the
complement to a link in $S^3$), 
then the elements 
$$A(G,\rho_\pi), \quad \D_{G,\ove\rho}~\in 
\Sigma^{-1}\L$$
are equal up to multiplication by 
an invertible element
of $ \Sigma^{-1}\L$.
\enpr
\Prf
As before we can assume that 
$\pi(g_1)=t\in R[t, t^{-1}]$.
Suppressing the first $n$ columns of the matrix
$\psi(\AA)$ we obtain an 
$n(s-1)\times n(s-1)$-matrix 
$\MM'$.
Up to invertible elements of 
$\Sigma^{-1}\L$ we have:
\bq\lb{f:alexx}
\det\MM'=\D_{G,\ove\r}.
\end{equation}

The  boundary operator
$\pr'_1$ in the localized
complex 
$$
\FF'_*=\Sigma^{-1}\FF_*
=\Sigma^{-1}\L\tens{\ZZZ G}
F_*
$$
is an epimorphism (see
\rrf{f:pfff} for the definition of $F_*$),
and it is easy to deduce that 
the second boundary operator 
$\pr_2'$ in this complex
is isomorphic to a homomorphism 
$(0,\phi):\L^{n(s-1)}\to
\L^n\oplus\L^{n(s-1)}$
where the matrix of $\phi$ equals $\MM'$.
Thus up to invertible elements of 
$\Sigma^{-1}\L$ we have:
\bq\lb{f:dett}
\det\MM'=I_{n(s-1)}(\pr_2) 
=J^{(2)}_1(\FF'_*)=J^{(2)}_1(\FF_*)=A(G,\r_\pi),
\end{equation}
and the proof is 
completed. $\qs$

\subsection{Relation with Fitting invariants: 
the multi-variable case}
\mlb{su:rel_det_multi}

In the case when $H\approx 
\ZZZ^k$ with $k\geq 2$
the ring $\L=R[H]$ is isomorphic to 
$R[t_1, t_1^{-1}, ... , t_k, t_k^{-1}]$.
The twisted Alexander polynomial 
in this case is an element
of the ring $\L$ itself 
(cf. \cite{Wada}, p.253).
Indeed, choose a presentation
$p=(g_1,\dots, g_s~|~h_1,\dots, h_l)$
of $G$ in such a way, that 
$\pi(g_1)=t_1, \pi(g_2)=t_2$.
Let 
$\r:G\to GL_R(V)$
be a right representation of the group $G$
where $V$ is a free 
finitely generated $R$-module. Put 
$$
P_1 = \det (1-t_1\r(g_1)), \quad 
P_2 = \det (1-t_2\r(g_2)). 
$$
Then
\bega
P_1=1+a_1t_1+...+a_nt_1^n, 
\mxx{ where } a_n=\pm\det \r(g_1);\\
P_2=1+b_1t_2+...+b_nt_2^n, 
\mxx{ where } b_n=\pm\det \r(g_2).
\end{gather}
Let $\AA=\AA(G,p)$
denote the Alexander matrix
for the presentation $p$.
For any sequence 
$S=i_1<i_2< ... <i_{n(s-1)}$
of integers put 
$$
\a_1^S=\det\psi(\AA_1)^S,\quad
\a_2^S=\det\psi(\AA_2)^S.
$$

The property \rrf{f:equall} implies
\bq\lb{f:m1_m2}
P_2\a_1^S
=
P_1\a_2^S \mxx{ for every } S.
\end{equation}
The polynomials 
$P_1, P_2$ are relatively prime in $\L$, 
therefore $P_1~|~\a_1^S$
in $\L$, and 
$$\D_{G,\ove\r}=\frac {Q_1}{P_1} \in \L,
\mxx{ where }
Q_1= {\rm GCD}_S~ \a_1^S, 
$$
as we have claimed above.
Moreover, it is easy to deduce from
\rrf{f:m1_m2}
 the following :
$$ {\rm GCD}_{S}\big(\a_1^S,\a_2^S\big) 
~\Big|~ \frac {Q_1}{P_1},
$$
and we obtain
\bepr\mlb{p:wada_det_multi1}
The Fitting invariant 
$A(G,\r_\pi)$
divides 
the \talp~ $\D_{G,\ove\r}$
in the ring $\L=R[H]$. $\qs$
\enpr

Now we proceed to  an analog 
of Proposition \mrf{p:alex_det_links_one}
for the multi-variable case. Let us introduce
the corresponding localization.
\bede\mlb{d:extremities}
Let $\m:H\to\RRR$ be any non-trivial group \ho.
\belis\item
We say that 
an element $x\in \L$
has {\it $\mu$-monic ends} 
if
\begin{gather*}
x=x_0h_0+\sum_{i=1}^{a-1} x_i h_i+ x_ah_a,
\mxx{ with }x_i\in R, h_i\in H, \\
\mx{ where } x_0, x_a \in R^\bu 
\mxx{ and }
\mu(h_0)< \mu(h_i) <\mu(h_a) \quad\forall i: 0<i<a.
\end{gather*}
\item 
The  multiplicative subset
of all the elements of 
$\L$ with $\mu$-monic
ends will be denoted
$\Sigma_\m\sbs \L$.
\enlis\end{defi}

\bepr\mlb{p:alex_det_links_multi}
Assume that the group $G$ has a 
presentation with 
$s$ generators and $s-1$ relators.
Then the images of the elements
$A(G,\r_\pi)$
and 
$\D_{G,\ove\r}$
in the ring $\Sigma_\mu^{-1} \L$
are equal up to invertible elements of this ring.
\enpr
\Prf
Similar to the proof of 
Proposition 
\mrf{p:alex_det_links_one}. $\qs$

It is natural to ask whether 
the Fitting invariant and 
\talp~are  equal up to 
invertible elements of $\L$,
at least in the case of knot groups. 
I do not know if it is true.


\section{Three-dimensional manifolds}
\mlb{s:three}

In this section we study the particular case 
case of $\smo$ compact three-dimensional 
manifolds $M$ of zero 
Euler characteristic.
We prove Theorem \mrf{t:van_nov_3}
which gives a criterion 
for vanishing of the twisted
Novikov homology of $M$
in terms of the \talp~ of $\pi_1(M)$
or, equivalently,
in terms of the first Fitting 
invariant of $\pi_1(M)$.
The second main result  of this section is 
Theorem \mrf{t:nov_van_conical}.
This theorem describes for a given 
right representation
$\rho$ of $\pi_1(M)$
and for a given epimorphism
$\pi_1(M)\to H$
the structure of
all classes $\xi\in H^1(M,\RRR)$
such that the $\rho$-twisted Novikov 
homology $\wh H_*(M,\r_\pi,\xi)$
vanishes. It turns out 
that this set is an
open conical polyhedral subset 
of $H^1(M,\RRR)$.
In Subsection
\mrf{su:nov_thur}
we discuss the relation  of the 
twisted Novikov
homology and the Thurston norm;
we suggest a natural 
question about this
relation. 

\subsection{The twisted Novikov homology 
and the twisted Alexander polynomial
for 3-manifolds}
\mlb{su:nov_3_alex}

Let $M$ be a compact $\smo$
three-dimensional manifold with $\chi(M)=0$.
Write $G=\pi_1(M),~\gL=\ZZZ G$.
We begin with two lemmas which describe 
the homotopy type of the chain complex
of the universal covering of $M$.
\bele\mlb{l:closed_ma}
Let $M$ be a closed connected 3-manifold.
Then the cellular chain complex 
of its universal covering
is chain homotopy equivalent to the following one:
\bq\lb{f:complex}
C_*
=
\{0\lTo\gL\lTo^{\pr_1} \gL^l\lTo^{\pr_2} 
\gL^l\lTo^{\pr_3} \gL\lTo 0\}
\end{equation}
where:
\belis\item
The matrix of $\pr_1$
equals
$$
\begin{pmatrix}
1-g_1 \\1-g_2 \\ \vdots \\ 1-g_l
\end{pmatrix}
$$
and   the elements $g_1,\dots,g_l\in G$ 
generate the group $G$.
\item
The matrix of $\pr_3$
equals
$$\left(
\begin{matrix}
1-\ve_1h_1,& 1-\ve_2h_2,& \dots & ,1-\ve_lh_l
\end{matrix}
\right)
$$
where $h_i\in G$ 
and $\ve_i=1$ if the loop $h_i$ 
conserves the
orientation, $\ve_i=-1$ if 
$h_i$ 
reverses the
orientation.
The elements $h_1,\dots,h_l\in G$ 
generate the group $G$.
\item
The matrix of $\pr_2$ is the 
Alexander matrix associated to some 
presentation of the group $G$.
\enlis
\enle
\Prf The lemma follows immediately from
the existence of the Hegaard decomposition 
for closed three-dimensional manifolds. $\qs$

\bele\mlb{l:chain_comp_not_closed}
Let $M$ be a  connected compact 3-manifold
with non-empty boundary and $\chi(M)=0$.
Then the  cellular chain complex of 
the  universal covering $\wi M$
is chain homotopy equivalent to the
following one:
\bq\lb{f:complex_non_closed}
C_*
=
\{0\lTo\gL\lTo^{\pr_1} \gL^l\lTo^{\pr_2} 
\gL^{l-1}\lTo 0\}
\end{equation}
where the matrix of $\pr_1$
equals
$$
\begin{pmatrix}
1-g_1 \\1-g_2 \\ \vdots \\ 1-g_l
\end{pmatrix}
$$
and   the elements $g_1, \dots, g_l\in G$ 
generate the group $G$.
The matrix of $\pr_2$ is the 
Alexander matrix associated to some 
presentation of the group $G$.
\enle
\Prf 
The Morse theory guarantees 
the existence of a Morse function
$f:M\to \RRR$
such that 
$$f|\pr M=const=\max_{x\in M} f(x)$$
and the number $m_i(f)$
of critical points of
index $i$ satisfies the following:
$$
m_0(f)=1, ~m_1(f)=m_2(f), ~m_3(f)=0.
$$
The cellular decomposition 
corresponding to $f$, 
satisfies the requirements
of the lemma. $\qs$

We shall use these lemmas to study
the \talp~and the Fitting invariants of $M$.
Let $V$ be a free \fg~$R$-module
(where $R$ is a commutative factorial ring)
and $\r:G\to GL_R(V)$ be a right 
representation.
Let $\pi:G\to H$
be an epimorphism of $G$
onto a free abelian \fg~group $H$.
Let $\xi:H\to\RRR$
be any non-trivial \ho.
Denote $R[H]$ by $\L$,
and let 
$\Sigma_\xi$ be the multiplicative subset
of Laurent polynomials with $\xi$-monic ends 
(see Definition 
\mrf{d:extremities}).
Let 

$$
\wi C_*=
\Sigma^{-1}_\xi 
\Big( V^H\tens{\gL} C_*\Big),
$$
where $V^H$
is endowed with the structure of 
a right $\gL$-module determined 
by the representation
$\r_\pi:G\to GL(V^H)$. Put 
$\Lxxi=\Sigma^{-1}_\xi\L$, and denote
$\rk V$ by $n$. We have a natural isomorphism
$$
\wi C_1
\approx
\Lxxi^{n(l-1)}\oplus \Lxxi^n
$$
(where the first summand of the direct sum 
correspond to the elements 
$e_1,...,e_{l-1}$
of the $\gL$-basis of $C_1$, 
and the second  
corresponds to $e_l$).
The projection of $\wi C_1$
onto the direct summand
$\Lxxi^{n(l-1)}$
will be denoted by $p_1$.
Similarly, 
the module 
$\wi C_2$
is naturally isomorphic 
to $\Lxxi^{n(l-1)}$ in 
the case $\pr M\not=\ems$
and to the direct sum 
$\Lxxi^{n(l-1)}\oplus \Lxxi^n$
in the case of closed manifolds.
Let 
\bq\lb{f:dd}
\DD=p_1\circ \Big( \pr_2|\Lxxi^{n(l-1)} \Big) .
\end{equation}

\bepr\mlb{p:homotopy_type}
Assume that 
\bq
\lb{f:conditionnn1}
\xi(g_l)<0,
\end{equation}
and in the case $\pr_1(M)=\ems$
assume moreover that 
\bq
\lb{f:conditionnn2}
\xi(h_l)<0.
\end{equation}
Then the chain complex 
$\wi C_*$
is chain homotopy equivalent 
to a free $\Lxxi$-chain complex
\bq\lb{f:homotopy_type}
0\lTo \Lxxi^{n(l-1)}
\lTo^\DD \Lxxi^{n(l-1)}
 \lTo 0,
\end{equation}
concentrated in dimensions 1 and 2.
\enpr
\Prf
Let us do the case of closed manifolds, 
the case
$\pr M\not=0$ is similar.
Observe that 
the homomorphisms
$$1-\rho_\pi(g_l),~ 1-\rho_\pi(h_l):
\Lxxi^n\to \Lxxi^n$$
are invertible.
Therefore, using a standard basis 
change in the chain complex
$\wi C_*$ we can split off from $\wi C_*$
two trivial chain complexes
$$
A_*\approx T_*(1, \Lxxi^n),\quad 
B_*\approx T_*(2,\Lxxi^n),$$
in such a way, that the resulting chain 
complex
is isomorphic to
\rrf{f:homotopy_type}.
$\qs$

Now we can establish the relation between 
the \talp~and the Fitting invariant 
of a three-manifold.
\bepr\mlb{p:talp_det_three}
Let $G$ be the fundamental group
of a compact three-dimensional $\smo$ manifold
$M$ with $\chi(M)=0$.
Let $\r:G\to GL_R(V)$ be a 
right representation,
where
$V$ is a free \fg~$R$-module
over a factorial commutative ring $R$. 
Let $\pi:G\to H$
be an epimorphism, where $H$ is a 
free abelian \fg~group. 
Let $\xi:H\to\RRR$ be any non-trivial \ho.

Then the Fitting invariant $A(G,\r_\pi)$ 
and the \talp~$\D_{G,\ove\r}$
are equal up to 
multiplication by an 
invertible element of
$\Lxxi$.
\enpr
\Prf
We shall do the case of closed manifolds, 
and the case $\pr M\not=\ems$
is even simpler, and will be omitted.

The Fitting invariant of a chain complex
depends only on 
its homotopy type and does not change
when we localize the base ring. 
Therefore 
we can use the chain complex 
\rrf{f:homotopy_type}
for computation of the image of 
$A(G,\r_\pi)$
in the ring $\Lxxi$
(we assume that the conditions 
\rrf{f:conditionnn1},
\rrf{f:conditionnn2}
 hold, which is easy to arrange 
by a permutation of the 
elements of the basis).
Thus the image of 
$A(G,\r_\pi)$
in the ring $\Lxxi$ equals $\det\DD$
(up to invertible elements of this ring).

To compute the \talp~
we shall use the matrix 
of the boundary operator $\pr_2$
from \rrf{f:complex}.
The $nl\times nl$-matrix $\DD'=\psi(\pr_2)$
satisfies
$\HH\circ \DD'=0$,
where $\HH=\psi\big((1-h_1,...1-h_l)\big).$

Let $\DD''$ denote the 
$nl\times n(l-1)$-matrix 
obtained by suppressing the last $n$
columns of the 
matrix of $\pr_2$.
Using the invertibility over $\Lxxi$
of the matrix $1-\rho_\pi(h_l)$
and the condition
$\DD''\circ\HH=0$,
we deduce that the last $n$ 
rows of the matrix
$\DD''$ are linear combinations of
the first $n(l-1)$ rows.
This implies
that the GCD of the 
$n(l-1)\times n(l-1)$-minors 
of the matrix $\DD''$
equals to the determinant of $\DD$.
We obtain therefore
the following equality:
$$
\D_{G,\bar\r}
=
\frac {\det\DD}{\det(1-\r_\pi(g_l))}~,
$$
and the proof is over. $\qs$

\beth\mlb{t:van_nov_3}
Let $M$ be a connected compact $\smo$ 
manifold of dimension 3 with $\chi(M)=0$.
Let $G=\pi_1(M)$, and let $\r:G\to GL_R(V)$
be a right representation of the group $G$,
where $V$ is a \fg~free module
over a commutative factorial ring  $R$.
Let $\pi:G\to H$ be an epimorphism onto 
a free abelian
\fg~group, and $\xi:H\to\RRR$ be a 
non-zero homomorphism.
Then the three following conditions 
are equivalent:
\belis\item 
The twisted  Novikov homology
$\wh H_i(M,\r_\pi,\xi)$
vanishes for all $i$.
\item The first twisted  
Novikov homology module
$\wh H_1(M,\r_\pi,\xi)$
vanishes.
\item
The Fitting invariant $A(G,\r_\pi)$
is $\xi$-monic.
\item
The \talp~ $\D_{G,\ove\r}$
is $\xi$-monic.
\enlis
\enth
\Prf
Observe that the two last conditions 
are equivalent in view of Proposition
\mrf{p:talp_det_three}.
Therefore it suffices  to 
prove that the first three conditions are 
equivalent.
Proposition \mrf{p:homotopy_type}
implies that 
the twisted Novikov homology 
$\wh H_*(M,\r_\pi,\xi)$
is isomorphic to the homology of the
chain complex 
\bq\lb{f:homotopy_type_nov}
0\lTo \Lxi^{n(l-1)}
\lTo^{\id \otimes \DD} \Lxi^{n(l-1)} \lTo 0,
\end{equation}
concentrated in dimensions 1 and 2, 
and $\det \DD$ equals to 
the Fitting invariant of 
$M$ \wrt~ $\r,\xi$.
The  homology of the chain complex 
\rrf{f:homotopy_type_nov}
vanishes for every $i$ if and only if 
it vanishes for $i=1$, and both 
these conditions are equivalent
to the invertibility of $\det \DD$ 
in $\Lxi$. This 
last condition 
holds if and only if this determinant is
is invertible in
$\Lx$. We have seen during the proof
of Proposition \mrf{p:talp_det_three}
that the elements  $\det \DD$ and $A(G,\rho_\pi)$
are equal up to multiplication by 
invertible elements of $\Lx$. 
Therefore the twisted Novikov 
homology vanishes if and only if
the element 
$A(G,\rho_\pi)$ is  $\xi$-monic.
$\qs$

The previous theorem allows to 
describe the structure of the set
of all $\xi$ such that the Novikov homology 
$\wh H_*(M,\r_\pi,\xi)$ vanishes.

\bede\mlb{d:acycc}
Let $G=\pi_1(M)$, and let 
$\r:G\to GL_R(V)$
be a right representation 
of the group $G$.
Let $\pi:G\to H$ be an 
epimorphism onto a 
free abelian
\fg~group.
A non-zero cohomology class 
$\xi\in H^1(M,\RRR)$ will 
be called {\it $(\r,\pi)$-acyclic},
if the twisted Novikov homology
$\wh H_1(M,\r_\pi,\xi)$ vanishes.
When the homomorphism $\pi$
is clear from the context,
we shall say that $\xi$ is $\r$-acyclic.
\end{defi}

\beth\mlb{t:nov_van_conical}
For a given right representation 
$\r:G\to GL_R(V)$ and a given epimorphism
$\pi:G\to H$
the set of all $(\r,\pi)$-acyclic classes 
$\xi\in H^1(M,\RRR)$
is 
an open polyhedral conical subset of 
$H^1(M,\RRR)$.
%
If $R$ is a field, then  
the set of all $(\r,\pi)$-acyclic classes
is either empty, or equals 
$H^1(M,\RRR)\sm\{0\}$
or is the complement in $H^1(M,\RRR)$
to a finite union of integral hyperplanes.
\enth
\Prf
The Fitting invariant $A(G,\r_\pi)$
is an element of $R[H]$.
The group $H$ is isomorphic 
to the integral lattice
$\ZZZ^k\sbs \RRR^k\approx H\otimes \RRR$, and
$A(G,\r_\pi)$
is then identified with a Laurent polynomial
$\AA$
in variables $t_1,...,t_k$ with
coefficients in $R$.
If $\AA=0$, then no class $\xi$
is $\r$-acyclic.
If $\AA$ is a monomial, 
$\AA=\a\cdot h$,
where $\a\in R$, and $h\in H$, then
either 
\been\item
$\a$ is invertible, and in this case
all classes $\xi$ are $\r$-acyclic, or
\item
$\a$ is non-invertible, 
and in this case there are no 
$\r$-acyclic classes. 
\enen
Now let us consider the non-degenerate case,
when the 
Newton polytope $\PP$
of the polynomial
$\AA$ contains  more than one point. 
 For a homomorphism
$\xi:H\to\RRR$
the polynomial $\AA$ is $\xi$-monic
if and only if the polytope $\PP$ 
has a vertex
$v$ such that:
\belis\item the coefficient $a_v$ of $\AA$
corresponding to this vertex 
is an invertible element of $R$,
\item
for every other vertex $v'$ of $\PP$ we have
$\xi(v)>\xi(v')$.
\enlis

For a given vertex $v$ the set
$\G_v$
of all $\xi$ satisfying the 
conditions 1) and 2) above
is an open polyhedral cone. 
Indeed, for a pair of vertices 
$v,v'$ of $\PP$
put
\begin{gather*}
\G_{v,v'}
=
\{\xi\in H^1(M,\RRR)~|~\xi(v)=\xi(v')\},
\\
\G^+_{v,v'}
=
\{\xi\in H^1(M,\RRR)~|~\xi(v)>\xi(v')\},
\end{gather*}
so that $\G^+_{v,v'}$
is one of the two open 
half-spaces corresponding to $\G_{v,v'}$.
Then
$$\G_v=\bigcap_{v'} \G^+_{v,v'}.$$
The sets $\G_v$ are open polyhedral cones
which are disjoint for 
different $v$. 

Now let $R$ be a field.
The  cases when $\AA=0$
or $\AA$ is invertible are 
done as above. 
When the Newton polytope $\PP$
of $\AA$ contains more than one vertex,
and the set of all $(\r,\pi)$-acyclic 
classes is
the complement in $H^1(M,\RRR)$
to the union of all the hyperplanes 
$\G_{v,v'}$. 
$\qs$

\subsection{Detecting fibred links}
\mlb{su:detect_fibr}

In this section we give a necessary condition for a 
link in $S^3$ to be fibred.
Let us first recall the  definition
of a fibred link and related notions.

\bede\mlb{d:top_categ}
\beli
\item Let $V$ be a compact topological 
$n-1$-manifold 
with $\pr V\not=\ems$
and let $h:V\to V$ be a homeomorphism
which restricts to the identity on $\pr V$.
Forming a mapping torus $V_h$, and identifying 
$(x,t)\sim (x,t')$ 
for each $x\in \pr V, t,t'\in S^1$
we obtain a closed topological manifold.
This manifold is denoted $B(V,h)$.

A closed manifold $M$
is called {\it an open book decomposition}
if it is homeomorphic to 
$B(V,h)$ for some $V$ and $h$.
The images in $M$ of fibers $V\times \{t\}\ (t\in S^1)$
are called {\it pages } of the open book,
and the image of $\pr V\times \{t\}$
in $M$ is called {\it binding }.

\item
A $\smo$-embedding of the disjoint union of
several copies of an oriented circle
into $S^3$
is called {\it oriented  link}.

\item
An oriented   link $L$ is called {\it fibred }
if there is an orientation preserving homeomorphism
$$
\phi: S^3\to B(F^2,h)
$$
where $F^2$ is a compact oriented surface
and the restriction of $\phi$ to $L$
is a preserving orientation homeomorphism
onto the binding $\pr F$ of the open book.
\enli
\end{defi}

For an oriented link $L$
let $G=\pi_1(S^3\sm L)$.
Let $\eta:G\to \ZZZ$
be a \ho, which sends every positive meridian
of $L$ to $1\in\ZZZ$.
Form the corresponding completion
$\wh\gL_\eta$ of the group ring $\gL=\ZZZ G$,
and let $\CC_*(\wi{S^3\sm L})$
be the cellular chain complex of
the universal covering of the link complement.

If the link $L$ is fibred then 
the complement $S^3\sm L$
admits a fibration over a circle.
Although the manifold $S^3\sm L$
is not compact it turns out that 
the corresponding analog of Theorem 
\mrf{t:nov_complex} holds; this is the 
subject of the next 
proposition.
\bepr\mlb{c:nov_fibr}
If the link $L$ is fibred, then
$$
H_*\Big( \wh\gL_\eta    \tens{\gL}   
\CC_*(\wi{S^3\sm L})  \Big)
=0.
$$
\enpr
\Prf
It is clear that $S^3\sm L$
is homotopy equivalent to 
the mapping torus $F^2_h$ 
where $h$ is a 
homeomorphism. The space 
$F^2_h$ 
is homotopy equivalent to the mapping torus
$F^2_g$
where $g:F^2\to F^2$
is a cellular map
(see \cite{ran_tors2}, Proposition 6.1).
Thus our proposition follows from
the next theorem.
\beth\mlb{t:tor_van}
Let $X$ be a finite connected 
CW-complex. Let $g:X\to X$
be a cellular map,
and let $X_g$ be the 
mapping torus. Let 
$G=\pi_1(X_g)$
and let $\eta:G\to\ZZZ$
be the homomorphism induced 
by the projection
$f:X_g\to S^1$.
Put $\gL=\ZZZ G$ and let 
$\wh\gL_\eta $
be the corresponding Novikov 
completion.
Then 
$$
H_*\Big( \wh\gL_\eta    \tens{\gL}   
\CC_*(\wi{X_g})  \Big)
=0.
$$
\enth
\Prf
Let 
\begin{gather*}
H=\ker\eta,\ R=\ZZZ H\\
P=\{\l\in\gL~|~\supp\l\in \eta^{-1}(]-\infty, 0]) \\
\wh P=\{\l\in \wh\gL_\eta~|~\supp\l\in 
\eta^{-1}(]-\infty, 0])\}.
\end{gather*}
Pick any $t\in G$ such that $\eta(t)=-1$.
Then the ring $P$ is isomorphic to 
the twisted polynomial
ring $R_\t[t]$,
and the ring $\wh P$
is isomorphic to 
the  twisted power series 
ring $R_\t[[t]]$
(where $\t$ is the isomorphism of the 
ring $R$ defined by 
$x\mapsto txt^{-1}$.)
Denote by $Y\to X_g$
the infinite cyclic covering
induced by the natural projection
$f:X_g\to S^1$
from the universal covering 
$\RRR\to S^1$.
The function $f$ lifts to a 
continuous function
$F:Y\to\RRR$;
for each $k\in \ZZZ$
the space
$F^{-1}([k,k+1])$
is homeomorphic to
the mapping cylinder $Z_g$
of the map $g:X\to X$.
Let $Y_k=F^{-1}(]-\infty, k])$.
For every $m<k$
the inclusion $Y_m\sbs Y_k$
is a homotopy equivalence.
Let $\wi Y$ be the universal covering
for $Y$, which is also 
a universal covering
for $X_g$.
Let $\wi Y_k$ be 
the inverse image of $Y_k$
in $\wi Y$.


The cellular chain complex
$\CC_*(\wi Y_k)$
is a free \fg~$P$-module.
We have for every $k$
$$
\wh\gL_\eta\tens{\gL}
\CC_*(\wi{Y})
=
\wh\gL_\eta\tens{\wh P}
\Big(\wh P\tens{P} 
\CC_*(\wi Y_k )\Big)
$$
thus it remains to show that
the chain complex
$\wh P\tens{P} \CC_*(\wi Y_0)$
is acyclic.
Observe that since the modules 
$\CC_*(\wi  Y_k  )$
are free \fg~$P$-modules,
we have
$$
\wh P\tens{P} \CC_*(\wi X_0)
=
\liminv \ \CC_*(\wi X_0, \wi X_{-k}   )
$$
(where $k\in \NNN$).
Since the inclusion
$X_{-k} \rInto X_0$
is a homotopy equivalence,
all  chain complexes
$\CC_*(\wi X_0, \wi X_{-k})$
are acyclic.

Applying Theorem A.19 of 
\cite{massey}, Appendix
we obtain:
$$
H_*\big(
\liminv \ \CC_*
(\wi Y_0, 
\wi Y_{-k})
\big)
\approx
\liminv \ H_*(\CC_*
(\wi Y_0, 
\wi Y_{-k}))=0, 
$$
and this completes the proof. $\qs$

Now let us apply this result to the twisted Novikov 
homology of the complement to a fibered link.
Let $R$ be a commutative ring.
Denote by $\xi$ the inclusion 
$\ZZZ\rInto \RRR$, and put $\L=R[\ZZZ]$.
 The ring $\Lxi$ is then  identified 
with the ring  $R((t))$, and 
the ring $\Lx$ is identified with the ring
$$
R\langle t\rangle =
\left\{ \frac {P(t)}{t^n\cdot (1+tQ(t))}
~|~
n\in\ZZZ \mxx{ and } P, Q\in R[t] \right\}
\sbs R((t)).
$$
The  $\xi$-monic elements
of $\L$ are identified with  polynomials
of the form
$\m(1+tQ(t))$ where $Q(t)\in R[t]$, and $\m\in R$
is invertible.
Invertible 
elements of $\Lx=R\langle t \rangle$
are also 
called {\it monic}.
Let $\r:G\to GL_R(V)$
be any right representation.

\bepr\mlb{p:fibr_twis}
If $L$ is fibred, 
then the twisted Novikov homology 
$\wh H_*(S^3\sm L,\r_\pi,\xi)$
equals to zero. 
\enpr
\Prf 
The proposition is deduced from 
Theorem \mrf{t:tor_van}
in the same way as
Theorem
\mrf{t:nov_complex_twi}
 is deduced from 
Theorem
\mrf{t:nov_complex}.
 $\qs$

Observe that the space
$S^3\sm L$ is homotopy equivalent 
to a compact 3-manifold 
with boundary. 
Indeed, let $N$ be an open  
tubular \nei~of $L$ in $S^3$.
Then $S^3\sm L $ has the same homotopy
type as $S^3\sm N$.
Theorem \mrf{t:van_nov_3}
implies then the following result. 
\beco\mlb{c:vanish_knots}
Let $G$ denote the group $\pi_1(S^3\sm L)$.
If the link $L$ is fibred, then:
\been\item
The Fitting invariant 
$A(G,\r_\pi)\in R[t, t^{-1}]$ is monic.
\item 
The \talp~$\D_{G,\ove{\r}}\in 
R\langle t \rangle$
is monic.
\enen
\enco

This result is related to 
the theorem due to H.Goda, T.Kitano, and 
T.Morifuji \cite{gokimo}.
Their theorem says that if a knot $K$
is fibred, and 
$\r:\pi_1(S^3\sm K)\to SL(n,F)$
is a representation (where $F$ is a field), 
then the leading coefficient 
of the \talp~associated 
to $\r$, equals to 1.

Our theorem is valid  in more general 
setting: it allows  representations in $GL(n,R)$,
where $R$ is a factorial ring, 
and not only in $SL(n, F)$.
On the other hand, for the representations
in $SL(n, F)$ the  theorem 
of \cite{gokimo} gives much more information
since it guarantees 
that the leading coefficient
of the \talp\ equals 1, and  
Corollary \mrf{c:vanish_knots}
asserts only that this coefficient is
non-zero, that is, the \talp~does not vanish.

\subsection{Thurston cones 
and $\r$-acyclicity cones}
\mlb{su:nov_thur}

Let $M$ be a closed three-dimensional 
$\smo$ manifold. ¨Put $G=\pi_1(M), ~\gL=\ZZZ G$.
Let $\VV(M)\sbs H^1(M,\RRR)$ be the 
subset of all the cohomology classes, 
representable by 
closed 1-forms without zeros.
Let 
$$\VV_h(M)~=
\{\xi\in H^1(M,\RRR)~|~
H_*(\wh\gL_\xi\tens{\gL} \CC_*(\wi M))
=0\}.$$
For a given right 
representation 
$\r: \pi_1(M)\to GL(\ZZZ^n)$
let $\VV_{alg}(M,\rho)$ 
be the subset of all 
$(\rho, \pi)$-acyclic
cohomology classes $\xi$,
where $\pi$ is the projection
$\pi_1(M)\to H_1(M)/Tors$.
It follows from the results of 
the previous section, that 
\bq\lb{f:inclus}
\VV(M)\sbs \VV_h(M)\sbs 
\VV_{alg}(M)=
\bigcap_{\rho\in \RR}\VV_{alg}(M,\rho),
\end{equation}
where $\RR$ is the set of 
all right representations 
$\pi_1(M)\to GL(\ZZZ^n)$.
The Thurston theorem \cite{th}
implies that the set $\VV(M)$ is an open 
conical  polyhedral subset of $H^1(M,\RRR)$.
The set $\VV_{alg}(M,\rho)$
is also an open conical  polyhedral subset 
of $H^1(M,\RRR)$, as it follows from the 
results of the previous section.
\pa
{\bf Question.}
For which manifolds $M$ the
inclusions \rrf{f:inclus} are equalities?
\pa
For every $\rho$ the set 
$\VV_{alg}(M,\rho)$
is effectively computable 
from  the twisted Alexander 
polynomial or from
the first Fitting invariant.
Thus, investigating  the 
inclusions \rrf{f:inclus}
will give a considerable amount 
of information 
on the structure of the set 
$\VV(M)$.
We think that computer experiments can 
help here, and can  
clarify the problem
(see the paper \cite{GodaPaj}
for an example of application of 
the Kodama's program KNOTS to a similar
question).
In the rest of this section we shall show that 
the properties of the inclusion
\rrf{f:inclus}
are quite sensitive to the class of 
the representations which we consider.
Let $\RR_\FF$
be the set of 
all representations $G\to GL(\FFF^n)$
where $\FFF$ is a finite field.
We are going to show that 
there are closed manifolds $M$ with 
$$\VV(M) \not=  \bigcap_{\rho\in\RR_\FF}\VV_{alg}(M,\rho).$$
Let
$M$ be any 3-manifold with 
a following property $(L)$:

\belis \item[L1)]
There is an open 
subset $U\sbs  H^1(M,\RRR)$
such that every element of $U$ 
does not contain a nowhere vanishing 1-form.
\item[L2)]
There is  a non-zero class 
$\xi_0\in H^1(M,\RRR)$
which contains a nowhere vanishing 1-form.
\enlis

Existence of such manifolds 
is indicated in \cite{th}, p. 125-127.
If $M$ has property $(L)$, then 
for every $\r$
the set $\VV(M,\r)$
of all $\r$-acyclic classes 
is non-empty, and by 
Theorem \mrf{t:nov_van_conical}
this set 
equals $H^1(M,\RRR)\sm\{0\}$ 
or is the complement 
in  $H^1(M,\RRR)$  
to a finite union of integral 
hyperplanes
(the case $\VV(M,\r)=\ems$
is excluded by the property L2)).
Therefore  $\VV(M,\r)$ is open and dense 
in $H^1(M,\RRR)$  
and the intersection $I$ of all 
the sets 
$\VV_{alg}(M,\r)$
over all right 
 representations in $GL(\FFF^n)$
(where $\FFF$ is a finite field), 
is a residual subset.
Thus $I$  
 intersects  every open subet 
of $H^1(M,\RRR)$,
in particular $U\cap I\not=\emptyset$,
and therefore there exist cohomology
classes $\xi\in H^1(M,\RRR)$
which are $\r$-acyclic for every 
representation $\r$ in $GL(\FFF^n)$,
but not representable by a nowhere 
vanishing closed 1-form.


\end{document}


\begin{thebibliography}{99}
\label{refer}




\bibitem{bz}
\bz

\bibitem{CockroftSwan}
\CockroftSwan{\twodim}

\bibitem{cohen}
\cohen


\bibitem{crfox}
\crfox

\bibitem{eisenbud}
\eisenbud


\bibitem{farber}
\farber


\bibitem{g1}
\Goda{\Heeg}

\bibitem{g2}
\Goda{\handlenumber}

\bibitem{gokimo}
\gokimo

\bibitem{gomo}
\gomo

\bibitem{GodaPaj}
\GodaPaj{\twn}








\bibitem{kitano}
\kitano

\bibitem{latour}
\latour



\bibitem{Lin}
\Lin

\bibitem{massey}
\massey




\bibitem{Mcd}
\McDonald{\linalg}



\bibitem{Mcm}
\McMullen{\alexthur}


\bibitem{novidok}
\novidok

\bibitem{noviuspe}
\noviuspe

\bibitem{pasbor}
\pasbor

\bibitem{pazamet}
\pazamet

\bibitem{patou}
\patou


\bibitem{pasur}
\pasur

\bibitem{paasym}
\paasym

\bibitem{pajandran}
\pajandran


\bibitem{prw}
\prw

\bibitem{ran_tors2}
\ranicki{\torsiontwo}



\bibitem{rolfsen}
\rolfsen

\bibitem{schuetz1}
\schuetz{\oneparam}

\bibitem{schuetz2}
\schuetz{\orbs}

\bibitem{Thang_repre}
\Thang{\repre}

\bibitem{th}
\Thurston{\norm}


\bibitem{turaev_alex}
\turaev{\alex}

\bibitem{turaev_two_comple}
\turaev{\twocomple}




\bibitem{Wada}
\Wada

\bibitem{waldhausen}
\waldhausen

\bibitem{witten}
\witten

\end{thebibliography}
\end{document}